\documentclass[12pt,a4paper]{book}
\pdfoutput=1

\usepackage{latexsym,fullpage,amsfonts,amssymb,amsmath,amscd,graphics,epic,mathrsfs}
\usepackage{amsthm}

\usepackage[all]{xy}
\usepackage{hyperref}

\usepackage[final]{pdfpages}

\newtheorem*{theorem*}{Theorem}
\newtheorem*{remark*}{Remark}
\newtheorem{lemma}{Lemma}[section]
\newtheorem{property}[lemma]{Property}
\newtheorem{proposition}[lemma]{Proposition}
\newtheorem{remark}[lemma]{Remark}
\newtheorem{example}[lemma]{Example}
\newtheorem{theorem}[lemma]{Theorem}
\newtheorem{definition}[lemma]{Definition}
\newtheorem{notation}[lemma]{Notation}
\newtheorem{corollary}[lemma]{Corollary}
\newtheorem{introtheorem}{Theorem}

\newtheorem{introcorollary}[introtheorem]{Corollary}
\newtheorem{introremark}[introtheorem]{Remark}

\begin{document}
\newcommand{\pperp}{\hbox{$\perp\hskip-6pt\perp$}}
\newcommand{\Sym}{\operatorname{Sym}}
\newcommand{\N}{{\mathbb N}}
\newcommand{\PP}{{\mathbb P}}
\newcommand{\Z}{{\mathbb Z}}
\newcommand{\Q}{{\mathbb Q}}
\newcommand{\R}{{\mathbb R}}
\newcommand{\C}{{\mathbb C}}
\newcommand{\K}{{\mathbb K}}
\newcommand{\F}{{\mathbb F}}
\newcommand{\eps}{{\varepsilon}}
\newcommand{\ko}{{\mathcal O}}
\newcommand{\wx}{{\widetilde x}}
\newcommand{\wz}{{\widetilde z}}
\newcommand{\wa}{{\widetilde a}}
\newcommand{\wF}{{\widetilde F}}
\newcommand{\wpsi}{{\widetilde \psi}}
\newcommand{\wG}{{\widetilde G}}
\newcommand{\wv}{{\widetilde v}}
\newcommand{\bz}{{\boldsymbol z}}
\newcommand{\bp}{{\boldsymbol p}}
\newcommand{\wy}{{\widetilde y}}
\newcommand{\we}{{\widetilde e}}
\newcommand{\wg}{{\widetilde g}}
\newcommand{\wb}{{\widetilde b}}
\newcommand{\wu}{{\widetilde u}}
\newcommand{\wc}{{\widetilde c}}
\newcommand{\bi}{{\omega}}
\newcommand{\bx}{{\boldsymbol x}}
\newcommand{\Log}{{\operatorname{Log}}}
\newcommand{\pr}{{\operatorname{pr}}}
\newcommand{\Graph}{{\operatorname{Graph}}}
\newcommand{\jet}{{\operatorname{jet}}}
\newcommand{\Tor}{{\operatorname{Tor}}}
\newcommand{\sqh}{{\operatorname{sqh}}}
\newcommand{\const}{{\operatorname{const}}}
\newcommand{\Arc}{{\operatorname{Arc}}}
\newcommand{\Sing}{{\operatorname{Sing}}}
\newcommand{\Span}{{\operatorname{Span}}}
\newcommand{\Aut}{{\operatorname{Aut}}}
\newcommand{\Int}{{\operatorname{Int}}}
\newcommand{\Aff}{{\operatorname{Aff}}}
\newcommand{\Area}{{\operatorname{Area}}}
\newcommand{\val}{{\operatorname{Val}}}
\newcommand{\conv}{{\operatorname{conv}}}
\newcommand{\rk}{{\operatorname{rk}}}
\newcommand{\ow}{{\overline w}}
\newcommand{\ov}{{\overline v}}
\newcommand{\ks}{{\cal S}}
\newcommand{\red}{{\operatorname{red}}}
\newcommand{\kc}{{\cal C}}
\newcommand{\ki}{{\cal I}}
\newcommand{\kj}{{\cal J}}
\newcommand{\ke}{{\cal E}}
\newcommand{\kz}{{\cal Z}}
\newcommand{\tet}{{\theta}}
\newcommand{\Del}{{\Delta}}
\newcommand{\bet}{{\beta}}
\newcommand{\mm}{{\mathfrak m}}
\newcommand{\kap}{{\kappa}}
\newcommand{\del}{{\delta}}
\newcommand{\sig}{{\sigma}}
\newcommand{\alp}{{\alpha}}
\newcommand{\Sig}{{\Sigma}}
\newcommand{\Gam}{{\Gamma}}
\newcommand{\gam}{{\gamma}}
\newcommand{\Lam}{{\Lambda}}
\newcommand{\lam}{{\lambda}}
\newcommand{\om}{{\omega}}

\newcommand{\ram}{\operatorname{ram}}
\newcommand{\tr}{\operatorname{Tr}}
\newcommand{\www}{\mathfrak{W}}
\newcommand{\ob}{\operatorname{ ob}}
\newcommand{\ind}{\operatorname{ind}}
\newcommand{\bc}{\operatorname{bc}}
\newcommand{\symm}{\mathbb{S}}
\newcommand{\res}{\operatorname{Res}_{s=1}}
\newcommand{\ad}{\operatorname{ad}}
\newcommand{\bes}{\mathfrak{B}}
\newcommand{\eisen}{\mathcal{E}}
\newcommand{\rest}{\lvert}
\newcommand{\dd}{\mathfrak{D}}
\newcommand{\uu}{\mathcal{U}}
\newcommand{\hh}{\mathcal{H}}
\newcommand{\ff}{\mathcal{F}}
\newcommand{\jj}{\mathfrak{J}}
\newcommand{\pp}{\mathfrak{p}}
\newcommand{\oo}{\mathcal{O}}
\renewcommand{\P}{\mathcal{P}}
\newcommand{\ch}{{\bf 1}}
\newcommand{\sss}{\mathcal{S}}
\newcommand{\ww}{\mathcal{W}}
\newcommand{\bb}{\mathcal{B}}
\newcommand{\ii}{\mathcal{I}}
\newcommand{\LL}{\mathcal{L}}
\newcommand{\Ker}{\mathfrak{K}}
\newcommand{\ccc}{\mathcal{C}}
\newcommand{\nm}{\operatorname{Nm}}
\newcommand{\bs}{\backslash}
\newcommand{\swrtz}{\mathcal{S}}
\newcommand{\Hom}{\operatorname{Hom}}
\newcommand{\diag}{\operatorname{diag}}
\newcommand{\Spec}{\operatorname{Spec}}
\newcommand{\Lie}{\operatorname{Lie}}

\newcommand{\sgn}{\operatorname{sgn}}
\newcommand{\kk}{{\bf K}}
\newcommand{\zz}{\mathbb{Z}}
\newcommand{\nn}{\mathbb{N}}
\newcommand{\cc}{\mathbb{C}}
\newcommand{\qq}{\mathbb{Q}}
\newcommand{\rr}{\mathbb{R}}
\newcommand{\disc}{\Delta}
\newcommand{\modulus}{\delta}
\newcommand{\aaa}{\mathfrak{a}}
\newcommand{\aut}{\mathcal{A}}
\newcommand{\A}{\mathbb{A}}
\newcommand{\I}{\mathbb{I}}
\newcommand{\vol}{\operatorname{vol}}
\newcommand{\Ht}{H}
\newcommand{\diff}{\mathfrak{d}}
\newcommand{\triv}{{\bf 1}}
\newcommand{\stab}{\operatorname{Stab}}
\newcommand{\Ind}{\operatorname{Ind}}
\newcommand{\id}{\operatorname{Id}}
\newcommand{\re}{\operatorname{Re}}
\renewcommand{\Im}{\operatorname{Im}}
\newcommand{\Fr}{\operatorname{Fr}}

\newcommand{\Gr}{\operatorname{Gr}}
\newcommand{\vrel}{{v_{\operatorname{rel}}}}
\newcommand{\vaux}{{v_{\operatorname{aux}}}}
\newcommand{\vcusp}{{v_{\operatorname{cusp}}}}

\newcommand{\match}[1]{\mathop{\longleftrightarrow}\limits^{#1}}
\newcommand{\E}{{\mathcal E}}
\newcommand{\Sc}{{\mathcal S}}
\newcommand{\G}{{\mathcal G}}
\newcommand{\T}{{\mathcal T}}
\newcommand{\Fre}{{Fr\'{e}chet \,}}
\newcommand{\et}{{\`{e}tale }}
\newcommand{\tG}{{\widetilde{G}}}
\newcommand{\tH}{{\widetilde{H}}}
\newcommand{\tO}{{\widetilde{O}}}
\newcommand{\Fou}{{\mathcal{F}}}
\newcommand{\abs}[1]{\left|{#1}\right|}
\newcommand{\norm}[1]{\lVert#1\rVert}
\newcommand{\sprod}[2]{\left\langle#1,#2\right\rangle}
\newcommand{\cC}{{C_c^{\infty}}}
\newcommand{\cD}{{\mathcal{D}}}
\newcommand{\g}{{\mathfrak{g}}}
\newcommand{\h}{{\mathfrak{h}}}
\newcommand{\z}{{\mathfrak{z}}}
\newcommand{\cgot}{{\mathfrak{c}}}
\newcommand{\Supp}{\mathrm{Supp}}

\newcommand{\uG}{{\underline{G}}}
\newcommand{\uH}{{\underline{H}}}
\newcommand{\uX}{{\underline{X}}}
\newcommand{\uU}{{\underline{U}}}
\newcommand{\uV}{{\underline{V}}}
\newcommand{\uW}{{\underline{W}}}
\newcommand{\uY}{{\underline{Y}}}
\newcommand{\uZ}{{\underline{Z}}}
\newcommand{\gd}{\g^{\sigma}}
\newcommand{\oF}{{\overline{F}}}
\newcommand{\cO}{{\mathcal{O}}}
\newcommand{\cN}{{\mathcal{N}}}
\newcommand{\Zar}{{Zariski }}
\renewcommand{\sl}{{\mathrm{sl} }}
\newcommand{\sll}{{\mathfrak{sl}}}
\newcommand{\Mp}{\operatorname{Mp}}
\newcommand{\GL}{\operatorname{GL}}
\newcommand{\charact}{\operatorname{char}}
\newcommand{\SL}{\operatorname{SL}}
\newcommand{\gl}{{\mathfrak{gl}}}
\newcommand{\HC}{\operatorname{HC}}

\newcommand{\rep}{{representation }}
\newcommand{\repp}{{representation. }}
\newcommand{\Rep}{{Representation }}
\newcommand{\reps}{{representations }}
\newcommand{\repsp}{{representations. }}
\newcommand{\irrep}{{irreducible representation }}
\newcommand{\irrepp}{{irreducible representation. }}
\newcommand{\irreps}{{irreducible representations }}
\newcommand{\irrepsp}{{irreducible representations. }}
\newcommand{\Irreps}{{Irreducible representations }}
\newcommand{\dist}{{distribution }}
\newcommand{\distp}{{distribution. }}
\newcommand{\dists}{{distributions }}
\newcommand{\distsp}{{distributions. }}
\newcommand{\Dist}{{Distribution }}
\newcommand{\Dists}{{Distributions }}

\newcommand{\oU}{\operatorname{U}}
\newcommand{\Sp}{{\mathrm{Sp}}}
\newcommand{\SO}{{\mathrm{SO}}}

%

%

%

\begin{titlepage}

\thispagestyle{empty}


\begin{center}
\large{Weizmann Institute of Science\\
Faculty of Mathematics and Computer Science}

\vspace{1.5cm} \Huge{\bf Multiplicity One Theorems \\ and
Invariant Distributions}

\vspace{1cm} \large{Thesis submitted for the degree\\
``Doctor of Philosophy"\\
by\\}

\vspace{1.5cm} {\Large \bf Dmitry Gourevitch }

\vspace{1.5cm} \large{Submitted to the Scientific Council of
the\\
Weizmann Institute of Science\\
Rehovot, Israel}

\vspace{1.5cm}

\large{Prepared at the Mathematics Department under the supervision of\\
Professor Joseph Bernstein and Professor Stephen Gelbart}

\vspace{1cm} {May 2009}

\end{center}
\end{titlepage}

\newpage
$ $\\\\\\\\ \thispagestyle{empty}
\begin{center}
{ \it \large To the memory of
my grandfathers\\
Naum Gourevitch and David Kagan}
\end{center}
\newpage
\begin{center}
$ $

\vspace{5cm}{\large This work was carried out under the
supervision of}

\vspace{2cm}{\Large \bf Professor Joseph Bernstein and\\ Professor
Stephen Gelbart}
\end{center}
\thispagestyle{empty}
\newpage
 $ $ \thispagestyle{empty}
\newpage
$ $ \thispagestyle{empty}

\section*{Acknowledgements}

First of all I would like to thank my advisors {\bf Joseph
Bernstein} and {\bf Stephen Gelbart} for guiding me through this
project.

I cordially thank {\bf Joseph Bernstein} for teaching me most of
the mathematics I know and for showing me his approach to
mathematics and to research in general.

I would like to thank my parents, {\bf Lena and Ilya} for teaching
me their approach to life and to work and for their support in my
carrier and in all decisions that I made in my life.

I am pleased to thank my wife {\bf Anna} for her advise, help and
support in life in general and in mathematics in particular. I
also cordially thank our children {\bf Mark} and {\bf Galina} for
their patience to me while I was busy working on this thesis.

I am glad to thank by best friend {\bf Avraham Aizenbud} with whom
large portion of this work was done in collaboration. We have been
studding and working together since 1997. His friendship always
supported me both in the darkest and in the happiest times in my
life.

I would like to thank my "scientific elder brother" {\bf Eitan
Sayag} for his encouragement and guiding in the last 13 years,
which is almost half of my life. I am very glad that my old wish
to collaborate with Eitan was fulfilled during my PhD studies.

I cordially thank {\bf Steve Rallis} and {\bf G\'{e}rard
Schiffmann} for releasing their ideas on the subject in their
preprint \cite{RS2}. This allowed us to write our preprint
\cite{AG_PfMult1Con}. My first paper on this subject, \cite{AGRS},
was a synthesis of these two preprints. I also thank \textbf{David
Soudry} for showing us the preprint \cite{RS2} just at the time we
needed it.

I would like to thank {\bf David Kazhdan} who drew my interest to
this subject by suggesting us to prove Theorem \ref{PGLDist} in
the Archimedean case.

During the years of my PhD studies I have benefitted from talks
with many people. I cannot list all of them, but I would like to
thank {\bf Moshe Baruch}, {\bf Vladimir Berkovich}, {\bf Patrick
Delorme}, \textbf{Gerrit van Dijk}, {\bf Yuval Flicker},
\textbf{Maria Gorelik}, \textbf{Nadia Gurevich}, \textbf{Herve
Jacquet}, {\bf David Kazhdan}, \textbf{Bernhard Kroetz}, {\bf Erez
Lapid}, {\bf Omer Offen}, {\bf Andrey Reznikov}, {\bf Siddhartha
Sahi}, {\bf Yiannis Sakellaridis}, {\bf Eugenii Shustin},
\textbf{David Soudry}, {\bf Ilya Tyomkin}, \textbf{Yakov
Varshavsky}, {\bf Nolan Wallach}, {\bf Wee Teck Gan} and
\textbf{Oksana Yakimova} for fruitful discussions.

I would also like to thank {\bf Semyon Alesker}, {\bf Erez Lapid}
and {\bf Binyong Sun} for useful remarks.

I wish to thank {\bf Vladimir Berkovich}, {\bf Maria Gorelik},
{\bf Dmitry Novikov} and {\bf Sergei Yakovenko} for their
encouragement, guidance and support.

Finally, I would like to thank the secretary of the pure math
department, {\bf Gizel Maymon}, and the previous secretary {\bf
Terry Debesh} for defending me from all kinds of burocracy.

During the work on this project I twice visited the Hausdorff
Center of Mathematics and Max Planck Institute fur Mathematik in
Bonn. An important part of this work was done during those visits.
I wish to thank both institutions for inspiring environment and
working conditions.

During the work on this thesis I was supported in part by a BSF
grant, a GIF grant, and an ISF center of excellency grant.


\newpage
\section*{Summary}
Let $F$ be a local field of characteristic zero, e.g. $F=\R$ or
$F=\Q_p$.

Consider the standard imbedding $\mathrm{GL}(n, F) \subset
\mathrm{GL}(n+1, F)$. Let $\mathrm{GL}(n, F)$ act on
$\mathrm{GL}(n+1, F)$ by conjugation.

The main goal of this thesis was to prove the following theorem.
\begin{theorem*} 
Every $\mathrm{GL}(n, F)$-invariant distribution on
$\mathrm{GL}(n+1, F)$ is invariant with respect to transposition.
\end{theorem*}

We proved it in \cite{AGRS} for non-Archimedean $F$ and in
\cite{AG_AMOT} for Archimedean $F$. For Archimedean $F$ it was
also proven independently, simultaneously and in a different way
in \cite{SZ}. The proof we present here is a combination of these
three approaches and works uniformly for both Archimedean and
non-Archimedean fields.

This theorem is important in representation theory, since it
implies the following multiplicity one theorem.
\begin{theorem*} 
Let $\pi$ be an irreducible admissible representation of
$\mathrm{GL}(n+1,F)$ and $\rho$ be an irreducible admissible
smooth representation of $\mathrm{GL}(n,F)$. Then
$$\dim \mathrm{Hom}_{\mathrm{GL}(n,F)}
(\pi|_{\mathrm{GL}(n,F)},\rho) \leq 1.$$
\end{theorem*}

Another corollary of the main theorem is the following one.

\begin{theorem*} 
Let $GL(n,F)$ act on itself by conjugation. Let $P(n,F) < GL(n,F)$
be the subgroup of matrices whose last row is $(0,\dots,0,1)$.
Then any $P(n,F)$-invariant distribution on $GL(n,F)$ is
$GL(n,F)$-invariant.
\end{theorem*}

This theorem in turn implies Kirillov's conjecture that states
that any irreducible unitary representation $GL(n,F)$ remains
irreducible unitary representation when restricted to $P(n,F)$. In
this way Kirillov's conjecture was originally proven in \cite{Ber}
(for non-archimedean $F$) and in \cite{Bar}(for archimedean $F$).

Analogs of our theorems hold also for the orthogonal groups.
Namely, let $V$ be a finite dimensional vector space over $F$ and
$q$ be any non-degenerate quadratic form on $V$. Extend $q$ to $V
\oplus F$ in the natural way and consider the embedding $O(V)
\hookrightarrow O(V \oplus F)$. Let $O(V)$ act on $O(V \oplus F)$
by conjugation.

\begin{theorem*} Every $O(V)$-invariant distribution on $O(V \oplus F)$ is
invariant with respect to transposition. Also, for any irreducible
admissible smooth representations $\pi$ of $O(V \oplus F)$ and
$\rho$ of $O(V)$ we have
$$\dim \mathrm{Hom}_{O(V)}
(\pi|_{O(V)},\rho) \leq 1.$$
\end{theorem*}
For non-archimedean $F$ this theorem proved in \cite{AGRS} and for
archimedean $F$ in \cite{SZ}. Both proofs use the theorem for
$GL(n,F)$.

In our proofs for non-archimedean $F$ we used tools dealing with
invariant distributions developed in \cite{BZ}, \cite{Ber} and
\cite{JR}. For archimedean $F$ we had to develop archimedean
analogs of some of these tools (see chapter \ref{Tools}).

However, the archimedean case is much more difficult and we needed
another, crucial, tool to solve it. This tool was the theory of
D-modules. In particular, we used a theorem on the singular
support of a D-module, see section \ref{SingSupp}. To make our
proof be uniform for both kinds of local fields we use the paper
\cite{Aiz} in which a non-Archimedean analog of this tool is
developed.


 \tableofcontents

 \setcounter{chapter}{-1}
\chapter{Introduction}

One may divide the task of \rep theory into two parts:

\begin{itemize}
\item Classify all \irreps of a given group
$G$.
\item Understand how a given \rep of $G$ "decomposes" into
irreducible ones.
\end{itemize}

For a \rep $\pi$ and an \irrep $\tau$ one can ask with what
multiplicity $\tau$ occurs in the decomposition of $\pi$. It is
interesting to know what happens if we consider an \irrep of $G$
as a representation of a subgroup $H \subset G$. Consider its
"decomposition" to irreducible \reps of $H$. It is said that the
pair $(G,H)$ has {\it multiplicity one} property if all those
\irreps appear in the decomposition with multiplicity one. Such
pairs are also called {\it strong Gelfand pairs}. If this property
holds just for the trivial representation instead of $\tau$ then
the pair $(G,H)$ is called a {\it Gelfand pair}.

Those notions were introduced by I.M. Gelfand in the 50s for pairs
of compact topological groups. In the 70s those notions were
developed by Gelfand and Kazhdan (in \cite{GK}) in the realm of
 reductive algebraic groups over non-Archimedean local fields (like the field of $p$-adic numbers). We recall their definitions and techniques,
 as well as an adaptation to Archimedean local fields (i.e. $\R$ and $\C$) in chapter \ref{ChGK}.

Let $F$ be a non-Archimedean local field of characteristic 0.
Consider the standard imbedding $\mathrm{GL}_n(F) \hookrightarrow
\mathrm{GL}_{n+1}(F)$. We consider the action of
$\mathrm{GL}_n(F)$ on $\mathrm{GL}_{n+1}(F)$ by conjugation.

In this thesis we prove
\begin{introtheorem}\label{StrongRep}
The pair $(\mathrm{GL}_{n+1}(F),\mathrm{GL}_n(F))$ is a strong
Gelfand pair.
\end{introtheorem}

Gelfand and Kazhdan also provided criteria for proving both
Gelfand properties. The criterion for strong Gelfand property is
the following. Suppose that we are given an anti-involution
$\sigma$ of $G$ that preserves $H$ and preserves all distributions
on $G$ that are invariant with respect to cojugation by elements
of $H$. Then $(G,H)$ is a strong Gelfand pair. By this criterion
Theorem \ref{StrongRep} follows from the following one:

\begin{introtheorem}\label{StrongDist}
Any $\mathrm{GL}_n(F)$ - invariant distribution on
$\mathrm{GL}_{n+1}(F)$ is invariant with respect to transposition.
\end{introtheorem}

For me, this theorem and invariant distributions in general are of
great interest on their own. Also, they have other applications in
representation theory. For example, Theorem \ref{StrongDist}
implies the following one.
\begin{introtheorem}\label{PGLDist}
Let $P_n \subset \mathrm{GL}_n$ be the subgroup consisting of all
matrices whose last row is $(0,...,0,1)$. Let $ \mathrm{GL}_n$ act
on itself by conjugation. Then every $ P_n(F)$ - invariant
distribution on $ \mathrm{GL}_n(F)$ is $ \mathrm{GL}_n(F)$ -
invariant.
\end{introtheorem}
\begin{proof}
Since $\GL_{n-1} < P_n$, any $ P_n(F)$ - invariant distribution
$\xi$ on $\GL_n(F)$ is  $\GL_{n-1}(F)$ - invariant. Hence by
Theorem \ref{StrongDist} $\xi$ is transposition invariant. Since
$P_n$ and $P_n^t$ generate $GL_n$, $\xi$ is $ \mathrm{GL}_n(F)$ -
invariant.
\end{proof}

This theorem in turn implies Kirillov's conjecture that states
that any irreducible unitary representation $GL(n, F)$ remains
irreducible unitary representation when restricted to $P_n(F)$. In
this way Kirillov's conjecture was originally proven by Bernstein
in \cite{Ber} for non-Archimedean $F$ and by Baruch in \cite{Bar}
for Archimedean $F$.

\section{The techniques and the structure of our proof}

Over non-Archimedean fields there is a powerful tool for proving
theorems on invariant distributions, due to Bernstein, Gelfand,
Kazhdan and Zelevinski.
  It says that if an algebraic group $G$ acts on an
algebraic variety $X$ and $\sigma$ is an involutive automorphism
of $X(F)$ that normalizes the action of $G$ and preserves all
$G(F)$-orbits then all $G(F)$-invariant distributions on $X(F)$
are $\sigma$-invariant. For a precise formulation see Theorem
\ref{Constructibility}.

This theorem is very powerful since it reduces a statement on
invariant distributions to a check of simple geometric conditions.
Many theorems on invariant distributions over non-Archimedean
fields were proven using this theorem.

However, in more difficult cases, including the ones under
consideration in this thesis, $\sigma$ preserves a majority of
$G(F)$-orbits but not all of them. In such cases there actually
exist $G(F)$-invariant but not $\sigma$-invariant distributions
defined on some locally closed subsets of $X(F)$ and one has to
prove that they cannot be invariantly continued to $X(F)$. For
this one has to use some non-geometric tools, like the ones
described in sections \ref{SecFour} and \ref{SingSupp}.

One can formulate a geometric condition which is necessary (though
not sufficient) in all cases. Namely, in order to preserve $G(F)$
invariant distributions $\sigma$ must preserve all {\bf closed}
$G(F)$-orbits.

In \cite{AG_Gen_HC_RJR} we formulated a sufficient condition that
one has to check for all closed orbits. This condition is not
geometric - it is a statement on invariant distributions on the
normal space to the orbit (see section \ref{SecDescent}). However,
the normal spaces to most orbits have dimension smaller than $X$
so one can hope that the problem on invariant distribution on them
is simpler. Also, the normal spaces are linear spaces which
enables to apply additional techniques, like the Fourier
transform. We called this method Generalized Harish-Chandra
descent, since in the case of a reductive group acting on itself
by conjugation it was developed by Harish-Chandra.

In the case of $GL_n$ acting on $GL_{n+1}$ by conjugation, for all
closed orbits except one, the problem we get on the normal space
is a product of problems of the same kind, possibly for other
local fields but with smaller $n$. Hence when proving the main
theorem by induction we may assume those problems to be solved,
and there is only one problematic orbit. This orbit consists of
one element - the identity. The normal space to it is the Lie
algebra $\gl_{n+1}(F)$. Using the fact that the problem is solved
near other closed orbits, we show that  any distribution $\xi$
which is $GL_n$-invariant and transposition anti-invariant must be
supported in the nilpotent cone. Since the Fourier transform of
$\xi$ has the same invariance properties, it also must be
supported in the nilpotent cone.

Now, when we have restrictions on the support of both $\xi$ and
its Fourier transform we can apply two kinds of "uncertainty
principles", that we develop in sections \ref{SecFour} and
\ref{SingSupp}. We apply them and show that such $\xi$ must be
zero. This is the most complicated part of the proof.

\section{Remarks}

Theorem \ref{StrongDist} is not proven yet for local fields of
positive characteristic (like $\F_q(t)$). The main difficulty is
that the Harish - Chandra descent technique does not work in this
case. We managed to prove a partial analog of Theorem
\ref{StrongDist}, see chapter \ref{ChWeak}.

Some chapters of this thesis are written uniformly for all local
fields. In other chapters or sections we assume that the field has
characteristic 0 and then we specify that in the beginning of the
chapter or section.

The first chapter of the thesis is dedicated to tools for working
with invariant distributions. Some of these tools existed before.
Many tools existed only in the non-Archimedean case and our
contribution was to develop their Archimedean analogs. However, we
developed the Generalized Harish-Chandra descent tool for all
local fields, and the tool of singular support (see section
\ref{SingSupp}) was  adapted to this kind of problems first in the
Archimedean case (in \cite{AG_AMOT}) and then continued to the
non-Archimedean case by Aizenbud in \cite{Aiz}.

Theorems \ref{StrongRep} and \ref{StrongDist} were first proven
for the non-Archimedean case in \cite{AGRS}. For the Archimedean
case it was proven simultaneously, independently and in different
ways in \cite{AG_AMOT} and \cite{SZ}. The proof we present here is
a synthesis of those three proofs. It works uniformly in both
cases thanks to \cite{Aiz}.

The main tools of the proof in \cite{AGRS} are the Homogeneity
Theorem (see section \ref{SecFour}) and certain non-linear
automorphisms $\nu_{\lambda}$ (see subsection \ref{subHC}). They
are used in the current proof for convenience, but are not
essential any more. The reason that they became less important for
us is that we have used a very powerful tool - the singular
support (see section \ref{SingSupp}). In the Archimedean case this
tool uses a deep result from the theory of $D$-modules. However,
the Homogeneity Theorem and the automorphisms $\nu_{\lambda}$ are
still used in an essential way in the proofs of multiplicity one
theorem for orthogonal groups. I hope that they will continue to
be used in further problems of this kind.

\section{Applications}

\begin{introremark} Using the tools developed here, one can show that Theorem \ref{StrongDist} implies an
analogous theorem for the unitary groups. Theorem \ref{StrongDist}
will be used during the Harish-Chandra descent and the main part
of the proof will be similar to the one we describe in chapter
\ref{ChProofMain}, since $U_n$ is a form of $GL_n$.
\end{introremark}

Analogs of Theorems \ref{StrongRep} and \ref{StrongDist} hold also
for the orthogonal groups. Namely, let $V$ be a finite dimensional
vector space over a local field $F$ of characteristic zero and $q$
be any non-degenerate quadratic form on $V$ . Extend $q$ to
$V\oplus F$ in the natural way and consider the embedding $O(V )
\hookrightarrow O(V\oplus F)$. Let $O(V )$ act on $O(V\oplus F)$
by conjugation.

\begin{introtheorem}\label{StrongODist}
Any $O(V)$ - invariant distribution on $O(V\oplus F)$ is invariant
with respect to transposition.
\end{introtheorem}

\begin{introcorollary}
$(O(V\oplus F),O(V ))$ is a strong Gelfand pair.
\end{introcorollary}

This theorem was conjectured by Bernstein and Rallis in the 1980s.
It was proven in \cite{AGRS} for non-Achimedean $F$ and in
\cite{SZ} for Archimedean $F$. Both proofs use Theorem
\ref{StrongDist}.

These proofs use the same techniques as the proof we present here,
except the singular support, which is not used in \cite{AGRS} and
only a partial analog of it is used in \cite{SZ}.

Theorem \ref{StrongDist} and its analog for unitary groups are
used in \cite{AGRS} and \cite{SZ} during Harish-Chandra descent
since the centralizer of a semisimple orthogonal operator is a
product of orthogonal, unitary and general linear groups over some
finite field extensions.

\begin{remark*}\label{SpecGroups}
One can ask whether analogous theorems hold for special linear,
unitary and orthogonal groups.

The pairs $(SL_{n+1},SL_n)$ and $(SU_{n+1},SU_n)$ are not strong
Gelfand pairs already for $n=1$ since $SL_1$ and $SU_1$ are
trivial and neither $SL_2$ nor $SU_2$ are commutative.

It is not known whether $(SO(V \oplus F),SO(V))$ is a strong
Gelfand pair. Of course, not every $SO(V)$ orbit is preserved by
transposition already for $\dim V = 1$ (since $SO(V)$ is trivial
in this case). On the other hand one can use another involution,
$\sigma_2$, obtained by composition of transposition and
conjugation by an element of $O(V) - SO(V)$. However, one can show
that there exists a closed $SO_4$-orbit in $SO_5$ that is not
transposition invariant and another closed $SO_4$-orbit that is
not $\sigma_2$-invariant.
\end{remark*}

In \cite{Sun} Theorem \ref{StrongDist}  is used to prove (at least
in the non-Archimedean case) a theorem on distributions on the
semidirect product of a symplectic group with its Heisenberg
group, invariant with respect to conjugation by the symplectic
group. That theorem in turn implies a certain multiplicity one
theorem involving symplectic and metaplectic groups,
which in turn is used in \cite{GGP} to show uniqueness of
Fourier-Jacobi models.
%

All multiplicity one theorems mentioned above can be used to prove splitting of periods of automorphic forms.\\

Assuming multiplicity at most one, a more difficult question is to
find when it is one. Some partial results are known.

For the orthogonal group (in fact the special orthogonal group)
this question has been studied by B.~Gross and D.~Prasad
(\cite{GP, P2}) who formulated a precise conjecture. For an up to
date account see \cite{GR,GGP,Wald}.

Multiplicity one theorems have important applications to the
relative trace formula, to automorphic descent,  to local and
global liftings of automorphic representations, and to
determinations of L-functions. In particular, multiplicity at most
one is used as a hypothesis in the work \cite{GPSR} in the study
of automorphic L-functions on classical groups.

\section{Structure of the thesis}
In chapter \ref{Tools} we develop tools to work with invariant
distributions.

In chapter \ref{ChGK} we give the definitions of the notions of
Gelfand pair and strong Gelfand pair. We also formulate the
Gelfand-Kazhdan criteria for these properties and prove them in
the Archimedean case. In particular, we show that Theorem
\ref{StrongDist} implies Theorem \ref{StrongRep}.

In chapter \ref{ChWeak} we prove Gelfand property for the pairs
$(\GL_{n+1}(F),\GL_n(F))$ and
$(\mathrm{O}_{n+1}(F),\mathrm{O}_n(F))$. For fields of
characteristic zero the results of this chapter are much weaker
then the main results. I decided to include this chapter for two
reasons. First, the proofs in this chapter work over any local
field, while the analogs of the main results for fields of
positive characteristic are not known. Second, the proofs in this
chapter are shorter and clearly show how to use some of the tools.

In chapter \ref{ChProofMain} we prove the main result, Theorem
\ref{StrongDist}.

\chapter{Tools to work with invariant distributions}\label{Tools}
This chapter is based on
\cite{BZ,Ber,Bar,AGS1,AG_Gen_HC_RJR,AG_AMOT,Aiz}. We will give
more precise bibliographical note in each section.

\section{Definitions}
\subsection{Conventions} \label{Conv}
\begin{itemize}
\item
Henceforth we fix a local field $F$. All the algebraic varieties
and algebraic groups that we will consider will be defined over
$F$.  In some sections and chapters we will assume that the
characteristic of $F$ is zero, and we will say so in the beginning
of the section or chapter. For simplicity we always assume that
the characteristic of $F$ is different from 2.
\item For a group $G$ acting on a set $X$ we denote by $X^G$ the set of fixed points of $X$. Also, for an element $x \in X$
we denote by $G_x$ the stabilizer of $x$.
\item By a reductive group we mean a (non-necessarily connected) algebraic reductive group.
\item We consider an algebraic variety $X$ defined over $F$ as an algebraic
variety over $\oF$ together with action of the Galois group
$Gal(\oF /F)$. On $X$ we only consider  the Zariski topology. On
$X(F)$ we only consider  the analytic (Hausdorff) topology. We
treat finite-dimensional linear spaces defined over $F$ as
algebraic varieties.
\item The tangent space of a manifold (algebraic, analytic, etc.) $X$ at $x$
will be denoted by $T_xX$.
 \item Usually
we will use the letters $X, Y, Z, \Delta$ to denote algebraic
varieties and the letters $G, H$ to denote reductive groups. We
will usually use the letters $V, W,U,K,M,N,C,O,S,T$ to denote
analytic spaces (such as $F$-points of algebraic varieties) and
the letter $K$ to denote analytic groups. Also we will use the
letters $L, V, W$ to denote vector spaces of all kinds.
\item For an algebraic variety $X$ defined over $\R$ we denote by
$X_{\C}$ the natural algebraic variety defined over $\R$  such
that $X_{\C}(\R)=X(\C)$. Note that over $\C$, $X_{\C}$ is
isomorphic to $X \times X$.
\item An action of a Lie algebra $\g$ on a (smooth, algebraic, etc) manifold $M$ is a Lie algebra homomorphism from $\g$ to the Lie algebra of vector fields on $M$.
Note that an action of a (Lie, algebraic, etc) group on $M$
defines an action of its Lie algebra on $M$.

\item For a Lie algebra $\g$ acting on $M$, an element $\alpha \in \g$ and a point $x \in M$ we denote by $\alpha(x) \in T_xM$ the value at point $x$ of the vector field corresponding to $\alpha$.
We denote by $\g x \subset T_xM$ or by $\g (x)$ the image of the
map $\alpha \mapsto \alpha(x)$ and by $\g_x \subset \g$ its
kernel.
\end{itemize}

\subsection{Analytic manifolds}
In this paper we  consider distributions over $l$-spaces, smooth
manifolds and Nash manifolds. $l$-spaces are locally compact
totally disconnected topological spaces and Nash manifolds are
semi-algebraic smooth manifolds.

For basic facts on $l$-spaces and distributions over them we refer
the reader to \cite[\S 1]{BZ}.

For basic facts on Nash manifolds and Schwartz functions and
distributions over them see subsection \ref{AppSubFrob} and
\cite{AG1}. In this paper we  consider only separated Nash
manifolds.

We  now introduce notation and terminology which  allows a uniform
treatment of the Archimedean and the non-Archimedean cases.

We will use the notion of an analytic manifold over a local field
(see e.g. \cite[Part II, Chapter III]{Ser}). When we say
"\textbf{analytic manifold}" we always mean analytic manifold over
some local field. Note that an analytic manifold over a
non-Archimedean field is in particular an $l$-space and an
analytic manifold over an Archimedean field is in particular a
smooth manifold.

\begin{definition}
A \textbf{B-analytic manifold} is either an analytic manifold over
a non-Archimedean local field, or a Nash manifold.
\end{definition}

\begin{example} If $X$ is a smooth algebraic variety,
then $X(F)$ is a B-analytic manifold and $(T_xX)(F) = T_x(X(F)).$
\end{example}

\begin{notation}
Let $M$ be an analytic manifold and $S$ be an analytic
submanifold. We denote by $N_S^M:=(T_M|_Y)/T_S $ the
\textbf{normal bundle to $S$ in $M$}. The \textbf{conormal bundle}
is defined by $CN_S^M:=(N_S^M)^*$. Denote by $\Sym^k(CN_S^M)$ the
k-th symmetric power of the conormal bundle. For a point $y\in S$
we denote by $N_{S,y}^M$ the normal space to $S$ in $M$ at the
point $y$ and by $CN_{S,y}^M$ the conormal space.
\end{notation}


Now let us introduce the term \emph{"vector system"}. This term
allows to formulate statements in wider generality.

\begin{definition}
For an analytic manifold $M$ we define the notions of a
\textbf{vector system} and a \textbf{B-vector system} over it.

For a smooth manifold $M$,  a vector system over $M$ is a pair
$(E,B)$ where $B$ is a smooth locally trivial fibration over $M$
and $E$ is a smooth (finite-dimensional) vector bundle over $B$.

For a Nash manifold $M$, a B-vector system over $M$ is a pair
$(E,B)$ where $B$ is a Nash fibration over $M$ and $E$ is a Nash
(finite-dimensional) vector bundle over $B$.

For an $l$-space $M$, a vector system over $M$ (or a B-vector
system over $M$) is a sheaf of complex linear spaces.
\end{definition}

In particular, in the case where $M$ is a point, a vector system
over $M$ is either a $\C$-vector space if $F$ is non-Archimedean,
or a smooth manifold together with a vector bundle in the case
where $F$ is Archimedean. The simplest example of a vector system
over a manifold $M$ is given by the following.

\begin{definition}
Let $\mathcal{V}$ be a vector system over a point $pt$. Let $M$ be
an analytic manifold. A \textbf{constant vector system with fiber
$\mathcal{V}$} is the pullback of $\mathcal{V}$ with respect to
the map $M \to pt$. We denote it by $\mathcal{V}_M$.
\end{definition}

\subsection{Distributions} \label{PrelDist}
\begin{definition}
Let $M$ be an analytic manifold over $F$. We define
$C_c^{\infty}(M)$ in the following way.

If $F$ is non-Archimedean then $C_c^{\infty}(M)$ is the space of
locally constant compactly supported complex valued functions on
$M$. We do not consider any topology on $C_c^{\infty}(M)$.

If $F$ is Archimedean then $C_c^{\infty}(M)$ is the space of
smooth compactly supported complex valued functions on $M$,
endowed with the standard topology.

For any analytic manifold $M$, we define the space of
distributions $\cD(M)$ by $\cD(M):=C_c^{\infty}(M)^*$. We consider
the weak topology on it.
\end{definition}

\begin{definition}
Let $M$ be a $B$-analytic manifold. We define $\Sc(M)$ in the
following way.

If $M$ is an analytic manifold over non-Archimedean field,
$\Sc(M):=C_c^{\infty}(M)$.

If $M$ is a Nash manifold, $\Sc(M)$ is the space of Schwartz
functions on $M$, namely smooth functions which are rapidly
decreasing together with all their derivatives. See \cite{AG1} for
the precise definition. We consider $\Sc(M)$ as a \Fre space.

For any $B$-analytic manifold $M$, we define the space of
\textbf{Schwartz distributions} $\Sc^*(M)$ by
$\Sc^*(M):=\Sc(M)^*$. Clearly, $\Sc(M)^*$ is naturally embedded
into $\cD(M)$.
\end{definition}

\begin{remark}
Schwartz distributions have the following two advantages over
general distributions:\\
(i) For a Nash manifold $X$ and an open Nash submanifold $U\subset
X$, we have the following exact sequence
$$0 \to \Sc^*_X(X \setminus U)\to \Sc^*(X) \to \Sc^*(U)\to 0.$$
(ii) Fourier transform defines an isomorphism $\Fou:\Sc^*(\R^n)
\to \Sc^*(\R^n)$.
\end{remark}
For a short survey on Schwartz functions and distributions on Nash
manifolds see section \ref{AppSubFrob}, and for more information
see \cite{AG1}.

\begin{notation}
Let $M$ be an analytic manifold. For a distribution $\xi \in
\cD(M)$ we denote by $\Supp(\xi)$ the support of $\xi$.

For a closed subset $N \subset M$ we denote
$$\cD_M(N):= \{\xi \in \cD(M)|\Supp(\xi) \subset N\}.$$
More generally, for a locally closed subset $N \subset M$ we
denote
$$\cD_M(N):=\cD_{M\setminus (\overline{N} \setminus N)}(N).$$

Similarly if $M$ is a $B$-analytic manifold and $N$ is a locally
closed subset we define $\Sc^*_M(N)$ in a similar vein.
\footnote{In the Archimedean case, locally closed is considered
with respect to the restricted topology - see section
\ref{AppSubFrob}.}
\end{notation}

\begin{definition}
Let $M$ be an analytic manifold over $F$ and $\E$ be a vector
system over $M$. We define $C_c^{\infty}(M,\E)$ in the following
way.

If $F$ is non-Archimedean then $C_c^{\infty}(M,\E)$ is the space
of compactly supported sections of $\E$.

If $F$ is Archimedean and $\E = (E,B)$ where $B$ is a fibration
over $M$ and $E$ is a vector bundle over $B$, then
$C_c^{\infty}(M,\E)$ is the complexification of the space of
smooth compactly supported sections of $E$ over $B$.

If $\mathcal{V}$ is a vector system over a point then we denote
$C_c^{\infty}(M,\mathcal{V}) := C_c^{\infty}(M,\mathcal{V}_M)$.
\end{definition}

We define $\cD(M,\E)$, $\cD_M(N,\E)$, $\Sc(M,\E)$, $\Sc^*(M,\E)$
and $\Sc^*_M(N,\E)$ in the natural way.

\section{Schwartz distributions on Nash manifolds}
\label{AppSubFrob}
In this section we give a short survey on Schwartz functions and
distributions on Nash manifolds. We will also prove some
properties of $K$-equivariant Schwartz distributions. We work in
the notation of \cite{AG1}, where one can read about Nash
manifolds and Schwartz distributions over them. More detailed
references on Nash manifolds are \cite{BCR} and \cite{Shi}.

Nash manifolds are equipped with the \textbf{restricted topology},
in which open sets are open semi-algebraic sets. This is not a
topology in the usual sense of the word as infinite unions of open
sets are not necessarily open sets in the restricted topology.
However, finite unions of open sets are open
 and therefore in the restricted topology we consider only
finite covers. In particular, if $E \to M$ is a Nash vector bundle
it means that there exists a \underline{finite} open cover $U_i$
of $M$ such that $E|_{U_i}$ is trivial.

\begin{notation}
Let $M$ be a Nash manifold. We denote by $D_M$ the Nash bundle of
densities on $M$. It is the natural bundle whose smooth sections
are smooth measures. For the precise definition see e.g.
\cite{AG1}.
\end{notation}

An important property of Nash manifolds is
\begin{theorem}[Local triviality of Nash manifolds; \cite{Shi}, Theorem I.5.12  ] \label{loctriv}
Any Nash manifold can be covered by a finite number of open
submanifolds Nash diffeomorphic to $\R^n$.
\end{theorem}

\begin{definition}
Let $M$ be a Nash manifold. We denote by $\G(M):= \Sc^*(M,D_M)$
the \textbf{space of Schwartz generalized functions} on $M$.
Similarly, for a Nash bundle $E \to M$ we denote by $\G(M,E):=
\Sc^*(M,E^* \otimes D_M)$ the \textbf{space of Schwartz
generalized sections} of $E$.

In the same way, for any smooth manifold $M$ we denote by
$C^{-\infty}(M):= \cD(M,D_M)$ the \textbf{space of generalized
functions} on $M$ and for a smooth bundle $E \to M$ we denote by
$C^{-\infty}(M,E):= \cD(M,E^* \otimes D_M)$ the \textbf{space of
generalized sections} of $E$.
\end{definition}

Usual $L^1$-functions can be interpreted as Schwartz generalized
functions but not as Schwartz distributions. We will need several
properties of Schwartz functions from \cite{AG1}.

\begin{property}[\cite{AG1}, Theorem 4.1.3] \label{pClass}  $\Sc(\R ^n)$ = Classical
Schwartz functions on $\R ^n$.
\end{property}

\begin{property}[\cite{AG1}, Theorem 5.4.3] \label{pOpenSet}
Let $U \subset M$  be a (semi-algebraic) open subset, then
$$\Sc(U,E) \cong \{\phi \in \Sc(M,E)| \quad \phi \text{ is 0 on } M
\setminus U \text{ with all derivatives} \}.$$
\end{property}

\begin{property}[see \cite{AG1}, \S 5]\label{pCosheaf}
Let $M$ be a Nash manifold. Let $M = \bigcup U_i$ be a finite open
cover of $M$. Then a function $f$ on $M$ is a Schwartz function if
and only if it can be written as $f= \sum \limits _{i=1}^n f_i$
where $f_i \in \Sc(U_i)$ (extended by zero to $M$).

Moreover, there exists a smooth partition of unity $1 =\sum
\limits _{i=1}^n \lambda_i$ such that for any Schwartz function $f
\in \Sc(M)$ the function $\lambda_i f$ is a Schwartz function on
$U_i$ (extended by zero to $M$).
\end{property}

\begin{property}[see \cite{AG1}, \S 5]\label{pSheaf}
Let $M$ be a Nash manifold and $E$ be a Nash bundle over it. Let
$M = \bigcup U_i$ be a finite open cover of $M$. Let $\xi_i \in
\G(U_i,E)$ such that $\xi_i|_{U_j} = \xi_j|_{U_i}$. Then there
exists a unique $\xi \in \G(M,E)$ such that $\xi|_{U_i} = \xi_i$.
\end{property}

We will also use the following notation.
\begin{notation}
Let $M$ be a metric space and $x \in M$. We denote by $B(x,r)$ the
open ball with center $x$ and radius $r$.
\end{notation}
\subsection{Submersion principle}
\begin{theorem}[\cite{AG2}, Theorem 2.4.16] \label{SurSubSec}
Let $M$ and $N$ be Nash manifolds and $s:M \rightarrow N$ be a
surjective submersive Nash map. Then locally it has a Nash
section, i.e. there exists a finite open cover $N= \bigcup \limits
_{i=1}^k U_i$ such that $s$ has a Nash section on each $U_i$.
\end{theorem}

\begin{corollary} \label{EtLocIsNash}
An \et map $\phi:M \to N$ of Nash manifolds is locally an
isomorphism. That means that there exists a finite cover $M =
\bigcup U_i$ such that $\phi|_{U_i}$ is an isomorphism onto its
open image.
\end{corollary}

\begin{theorem} \label{NashEquivSub}
Let $p:M \to N$ be a Nash submersion of Nash manifolds. Then there
exist a finite  open (semi-algebraic) cover $M = \bigcup U_i$ and
isomorphisms $\phi_i:U_i \cong W_i$ and $\psi_i:p(U_i) \cong V_i$
where $W_i\subset \R^{d_i}$ and $V_i \subset \R^{k_i}$ are open
(semi-algebraic) subsets, $k_i \leq d_i$ and $p|_{U_i}$ correspond
to the standard projections.
\end{theorem}
\begin{proof}
The problem is local, hence without loss of generality we can
assume that $N = \R^k$, $M$ is an equidimensional closed
submanifold of $\R^n$ of dimension $d$, $d \geq k$, and $p$ is
given by the standard projection $\R^n \to \R^k$.

Let $\Omega$ be the set of all coordinate subspaces of $\R^n$ of
dimension $d$ which contain $N$. For any $V \in \Omega$ consider
the projection $pr:M \rightarrow V$. Define $U_V= \{x \in M | d_x
pr$ is an isomorphism $\}$. It is easy to see that $pr|_{U_V}$ is
\et and $\{U_V\}_{V \in \Omega}$ gives a finite cover of $M$. The
theorem now follows from the previous corollary (Corollary
\ref{EtLocIsNash}).
\end{proof}

\begin{theorem} \label{NashSub}
Let $\phi:M \to N$ be a Nash submersion of Nash manifolds. Let $E$
be a Nash bundle over $N$. Then\\
(i) there exists a unique continuous linear map
$\phi_*:\Sc(M,\phi^*(E)\otimes D_M) \to \Sc(N,E \otimes D_N)$ such
that for any $f \in \Sc(N,E^*)$ and $\mu \in \Sc(M,\phi^*(E)
\otimes D_M)$ we have $$\int_{x \in N} \langle f(x),\phi_*\mu(x)
\rangle = \int_{x \in M} \langle \phi^*f(x), \mu(x) \rangle.$$ In
particular, we mean that both integrals converge. \\
(ii) If $\phi$ is surjective then $\phi_*$ is surjective.
\end{theorem}
\begin{proof}
$ $\\ \indent (i)

Step 1. Proof for the case when $M= \R^n$, $N= \R^k$, $k \leq n$,
$\phi$ is the
standard projection and $E$ is trivial.\\
Fix Haar measure on $\R$ and identify $D_{\R^l}$ with the trivial
bundle for any $l$. Define $$\phi_*(f)(x):= \int _{y \in \R^{n-k}}
f(x,y)dy.$$ Convergence of the integral and the fact that
$\phi_*(f)$ is a Schwartz function follows from standard calculus.

Step 2. Proof for the case when  $M \subset \R^n$ and $N \subset
\R^k$ are open
(semi-algebraic) subsets, $\phi$ is the standard projection and $E$ is trivial.\\
Follows from the previous step and Property \ref{pOpenSet}.

Step 3. Proof for the case when  $E$ is trivial.\\
Follows from the previous step, Theorem \ref{NashEquivSub} and
partition of unity (Property \ref{pCosheaf}).

Step 4. Proof in the general case.\\
Follows from the previous step and partition of unity (Property
\ref{pCosheaf}).

(ii) The proof is the same as in (i) except of Step 2. Let us
prove (ii) in the case of Step 2. Again, fix Haar measure on $\R$
and identify $D_{\R^l}$ with the trivial bundle for any $l$. By
Theorem \ref{SurSubSec} and partition of unity (Property
\ref{pCosheaf}) we can assume that there exists a Nash section
$\nu:N \to M$. We can write $\nu$ in the form $\nu(x) = (x,s(x))$.

For any $x \in N$ define $R(x):= \sup\{r \in \R_{\geq 0}|
B(\nu(x),r) \subset M \}$. Clearly, $R$ is continuous and
positive. By Tarski - Seidenberg principle (see e.g. \cite[Theorem
2.2.3]{AG1}) it is semi-algebraic. Hence (by \cite[Lemma
A.2.1]{AG1}) there exists a positive Nash function $r(x)$ such
that $r(x) < R(x)$. Let $\rho \in \Sc(\R^{n-k})$ such that $\rho$
is supported in the unit ball and its integral is 1. Now let $f
\in \Sc(N)$. Let $g \in C^{\infty}(M)$ defined by $g(x,y):=
f(x)\rho((y-s(x))/r(x))/r(x)$ where $x \in N$ and $y \in
\R^{n-k}$. It is easy to see that $g \in \Sc(M)$ and $\phi_*g=f$.
\end{proof}

\begin{notation}
Let $\phi:M \to N$ be a Nash submersion of Nash manifolds. Let $E$
be a bundle on $N$. We denote by $\phi^*:\G(N,E) \to
\G(M,\phi^*(E))$ the dual map to $\phi_*$.
\end{notation}

\begin{remark}
Clearly, the map $\phi^*:\G(N,E) \to \G(M,\phi^*(E))$ extends to
the map $\phi^*:C^{-\infty}(N,E) \to C^{-\infty}(M,\phi^*(E))$
described in \cite[Theorem A.0.4]{AGS1}.
\end{remark}

\begin{proposition} \label{EnoughPull}
Let $\phi:M \to N$ be a surjective Nash submersion of Nash
manifolds. Let $E$ be a bundle on $N$. Let $\xi \in
C^{-\infty}(N)$. Suppose that $\phi^*(\xi) \in \G(M)$. Then $\xi
\in \G(N)$.
\end{proposition}
\begin{proof}
It follows from Theorem \ref{NashSub} and Banach Open Map Theorem
(see \cite[Theorem 2.11]{Rud}).
\end{proof}

\section{Basic tools}
This section is based on \cite{BZ,Ber,Bar,AGS1}.
\begin{theorem}\label{Filt_nonarch}
Let an $l$-group $K$ act on an $l$-space $M$. Let $M =
\bigcup_{i=0}^l M_i$ be a $K$-invariant stratification of $M$. Let
$\chi$ be a character of $K$. Suppose that
$\Sc^*(M_i)^{K,\chi}=0$. Then $\Sc^*(M)^{K,\chi}=0$.
\end{theorem}
This theorem is a direct corollary of \cite[Corollary 1.9]{BZ}.

\begin{theorem} \label{NashFilt}
Let a Nash group $K$ act on a Nash manifold $M$. Let $N$ be a
locally closed subset. Let $N = \bigcup_{i=0}^l N_i$ be a Nash
$K$-invariant stratification of $N$. Let $\chi$ be a character of
$K$. Suppose that for any $k \in \Z_{\geq 0}$ and $0 \leq i \leq
l$, $$\Sc^*(N_i,\Sym^k(CN_{N_i}^M))^{K,\chi}=0.$$ Then
$\Sc^*_M(N)^{K,\chi}=0.$
\end{theorem}
For the proof see e.g. \cite[\S B.2]{AGS1}.

The following proposition sometimes helps to verify the conditions
of this theorem.
\begin{proposition} \label{Trick}
Let a Nash group $K$ act on a Nash manifold $M$. Let $V$ be a real
finite dimensional representation of $K$. Suppose that $K$
preserves the Haar measure on $V$. Let $U\subset V$ be an open
non-empty $K$-invariant subset. Let $\chi$ be a character of $K$.
Suppose that $\Sc^*(M \times U)^{K,\chi}=0$. Then $\Sc^*(M,
Sym^k(V))^{K,\chi}=0$.
\end{proposition}
For proof see \cite[Section B.4]{AGS1}.

\begin{theorem}[Frobenius descent] \label{Frob}
Let an analytic group $K$ act on an analytic manifold $M$. Let $N$
be an analytic manifold with a transitive action of $K$. Let
$\phi:M \to N$ be a $K$-equivariant map.

Let $z \in N$ be a point and $M_z:= \phi^{-1}(z)$ be its fiber.
Let $K_z$ be the stabilizer of $z$ in $K$. Let $\Delta_K$ and
$\Delta_{K_z}$ be the modular characters of $K$ and $K_z$.

Let $\E$ be a $K$-equivariant vector system over $M$. Then\\
(i) there exists a canonical isomorphism $$\Fr: \cD(M_z,\E|_{M_z}
\otimes \Delta_K|_{K_z} \cdot \Delta_{K_z}^{-1})^{K_z} \cong
\cD(M,\E)^K.$$ In particular, $\Fr$ commutes with restrictions to
open sets.

(ii) For B-analytic manifolds $\Fr$ maps $\Sc^*(M_z,\E|_{M_z}
\otimes \Delta_K|_{K_z} \cdot \Delta_{K_z}^{-1})^{K_z}$ to
$\Sc^*(M,\E)^K$.
\end{theorem}

For the proof of (i) see \cite[\S\S 1.5]{Ber}  and \cite[\S\S 2.21
- 2.36]{BZ}  for the case of $l$-spaces and \cite[Theorem
4.2.3]{AGS1} or \cite{Bar} for smooth manifolds. For the proof of
(ii) see \cite{AG_Gen_HC_RJR}, Appendix B.

The following fundamental theorem is a combination of
\cite[Theorem 6.13 and Theorem 6.15]{BZ}.
\begin{theorem}\label{Constructibility}
Suppose that $F$ is non-Archimedean. Let $\gamma$ be an algebraic
action of a linear algebraic group $G$ on an algebraic variety
$X$. Let $\sigma:X(F) \to X(F)$ be a homemorphism such that:\\
(i) For each $g \in G(F)$ there exists $G^{\sigma} \in G(F)$ such
that $\gamma(g)\sigma = \sigma \gamma(g^{\sigma})$.\\
(ii) For some number $n$ and $g_0 \in G(F)$, $\sigma^n =
\gamma(g_0)$.\\
(iii) $\sigma$ preserves all $G(F)$ orbits.

Then each $G(F)$-invariant distribution on $X(F)$ is invariant
under $\sigma$.
\end{theorem}

\begin{remark}
This theorem is very powerful since it reduces a statement on
 invariant distributions to check of simple geometric conditions.
 Many theorems on invariant distributions over non-Archimedean
 fields were proven using this theorem.

 However, in more difficult
 cases, including the ones under consideration in this thesis,
 $\sigma$ preserves majority of $G(F)$-orbits but not all of them.
 In such cases there actually exist $G(F)$-invariant but not
 $\sigma$-invariant distributions defined on some locally closed
 subsets of $X(F)$ and one has to prove that they cannot be
 invariantly continued to $X(F)$. For this one has to use some
 non-geometric tools, like the ones described in sections
 \ref{SecFour} and \ref{SingSupp}.
\end{remark}

We will also use the following straightforward proposition.

\begin{proposition} \label{Product}
Let $K_i$ be analytic groups acting on analytic manifolds $M_i$
for $i=1 \ldots n$. Let $\Omega_i \subset K_i$ be analytic
subgroups. Let $\E_i \to M_i$ be $K_i$-equivariant vector systems.
Suppose that
$$\cD(M_i,E_i)^{\Omega_i}=\cD(M_i,E_i)^{K_i}$$ for all $i$. Then
$$\cD(\prod M_i, \boxtimes E_i)^{\prod \Omega_i}=\cD(\prod M_i,
\boxtimes E_i)^{\prod K_i},$$ where $\boxtimes$ denotes the
external product.

Moreover, if $\, \Omega_i$, $K_i$, $M_i$ and $\E_i$ are
$B$-analytic then the analogous statement holds for Schwartz
distributions.
\end{proposition}
For the proof see e.g. \cite[proof of Proposition 3.1.5]{AGS1}.

\section{Generalized Harish-Chandra descent} \label{SecDescent}


Harish-Chandra developed a technique based on Jordan decomposition
that allows to reduce certain statements on conjugation invariant
distributions on a reductive group to the set of unipotent
elements, provided that the statement is known for certain
subgroups (see e.g. \cite{HCh}).

In this section we generalize an aspect of this technique to the
setting of a reductive group acting on a smooth affine algebraic
variety, using the Luna Slice Theorem. Our technique is oriented
towards proving Gelfand property for pairs of reductive groups.

In this section we assume that the characteristic of the local
field $F$ is zero. This section is based on \cite{AG_Gen_HC_RJR}.

We start with some preliminaries from algebraic geometry and
invariant theory.

\subsubsection{Categorical quotient}
\begin{definition}
Let an algebraic group $G$ act on an algebraic variety $X$. A pair
consisting of an algebraic variety $Y$ and a $G$-invariant
morphism $\pi:X \to Y$ is called \textbf{the quotient of $X$ by
the action of $G$} if for any pair $(\pi ', Y')$, there exists a
unique morphism $\phi : Y \to Y'$ such that $\pi ' = \phi \circ
\pi$. Clearly, if such pair exists it is unique up to a canonical
isomorphism. We will denote it by $(\pi _X, X/G)  $.
\end{definition}

\begin{theorem}[cf. \cite{Dre}] \label{Quotient}
Let a reductive group $G$ act on an affine variety $X$. Then the
quotient $X/G$ exists, and every fiber of the quotient map $\pi_X$
contains a unique closed orbit. In fact, $X/G:=\Spec\cO(X)^G.$
\end{theorem}

\subsection{Preliminaries on algebraic geometry over local fields}
\label{PrelLoc}
\subsubsection{$G$-orbits on $X$ and $G(F)$-orbits on $X(F)$} 

\begin{lemma}[\cite{AG_Gen_HC_RJR}, Lemma 2.3.4] \label{OrbitIsOpen}
Let $G$ be an algebraic group and let $H \subset G$ be a closed
subgroup. Then $G(F)/H(F)$ is open and closed in $(G/H)(F)$.
\end{lemma}

\begin{corollary}
Let an algebraic group $G$ act on an algebraic variety $X$. Let $x
\in X(F)$. Then $$N_{Gx,x}^{X}(F) \cong N_{G(F)x,x}^{X(F)}.$$
\end{corollary}

\begin{proposition} \label{LocClosedOrbit}
Let an algebraic group $G$ act on an algebraic variety $X$.
Suppose that $S \subset X(F)$ is a non-empty closed
$G(F)$-invariant subset. Then $S$ contains a closed orbit.
\end{proposition}
\begin{proof}
The proof is by Noetherian induction on $X$. Choose $x \in S$.
Consider $Z:=\overline{Gx} -Gx$.

If $Z(F)\cap S$ is empty then $Gx(F) \cap S$ is closed and hence
$G(F)x \cap S$ is closed by Lemma \ref{OrbitIsOpen}. Therefore
$G(F)x$ is closed.

If $Z(F)\cap S$ is non-empty then $Z(F) \cap S$ contains a closed
orbit by the induction assumption.
\end{proof}

\begin{corollary} \label{OpenClosedAll}
Let an algebraic group $G$ act on an algebraic variety $X$. Let
$U$ be an open $G(F)$-invariant subset of $X(F)$. Suppose that $U$
contains all closed $G(F)$-orbits. Then $U=X(F)$.
\end{corollary}

\begin{theorem}[\cite{RR}, \S 2 fact A, pages 108-109] \label{LocZarClosed}%
Let a reductive group $G$ act on an affine variety $X$. Let $x
\in X(F)$. Then the following are equivalent:\\
(i) $G(F)x \subset X(F)$ is closed (in the analytic topology).\\
(ii) $Gx\subset X$ is closed (in the Zariski topology).
\end{theorem}

\begin{definition}
Let a reductive group $G$ act on an affine variety $X$. We call an
element $x \in X$ \textbf{$G$-semisimple} if its orbit $Gx$ is
closed.
\end{definition}

In particular, in the case where $G$ acts on itself by
conjugation, the notion of $G$-semisimplicity coincides with the
usual one.

\begin{notation}
Let $V$ be an $F$-rational finite-dimensional representation of a
reductive group $G$. We set
$$Q_G(V):=Q(V):=(V/V^G)(F).$$ Since $G$ is reductive, there is a
canonical embedding $Q(V) \hookrightarrow V(F)$. Let $\pi : V(F)
\to (V/G)(F)$ be the natural map. We set$$\Gamma_G(V):=
\Gamma(V):= \pi^{-1}(\pi(0)).$$ Note that $\Gamma(V) \subset
Q(V)$. We also set$$R_G(V):= R(V):= Q(V) - \Gamma(V).$$
\end{notation}

\begin{notation}
Let a reductive group $G$ act on an affine variety $X$. For a
$G$-semisimple element $x \in X(F)$ we set $$S_x := \{y \in X(F)
\, | \, \overline{G(F)y} \ni x\}.$$
\end{notation}

\begin{lemma} \label{Gamma}
Let $V$ be an  $F$-rational finite-dimensional representation of a
reductive group $G$. Then $\Gamma(V) = S_0$.
\end{lemma}
This lemma follows from \cite[fact A on page 108]{RR} for
non-Archimedean $F$ and \cite[Theorem 5.2 on page 459]{Brk} for
Archimedean $F$.

\begin{example}
Let a reductive group $G$ act on its Lie algebra $\g$ by the
adjoint action. Then $\Gamma(\g)$ is the set of nilpotent elements
of $\g$.
\end{example}

\begin{proposition}
Let a reductive group $G$ act on an affine variety $X$. Let $x,z
\in X(F)$ be $G$-semisimple elements which do not lie in the same
orbit of $G(F)$. Then there exist disjoint $G(F)$-invariant open
neighborhoods $U_x$ of $x$ and $U_z$ of $z$.
\end{proposition}
For the proof of this Proposition see \cite{Lun2} for Archimedean
$F$ and \cite[fact B on page 109]{RR} for non-Archimedean $F$.

\begin{corollary} \label{EquivClassClosed}
Let a reductive group $G$ act on an affine variety $X$. Suppose
that $x \in X(F)$ is a $G$-semisimple element. Then the set $S_x$
is closed.
\end{corollary}
\begin{proof}
Let $y \in \overline{S_x}$. By Proposition \ref{LocClosedOrbit},
$\overline{G(F)y}$ contains a closed orbit $G(F)z$. If $G(F)z =
G(F)x$ then $y \in S_x$. Otherwise, choose disjoint open
$G$-invariant neighborhoods $U_z$ of $z$ and $U_x$ of $x$. Since
$z \in \overline{G(F)y}$, $U_z$ intersects $G(F)y$ and hence
contains $y$. Since  $y \in \overline{S_x}$, this means that $U_z$
intersects $S_x$. Let $t \in U_z \cap S_x$. Since $U_z$ is
$G(F)$-invariant, $G(F)t \subset U_z$. By the definition of $S_x$,
$x \in \overline{G(F)t}$ and hence $x \in \overline{U_z}$. Hence
$U_z$ intersects $U_x$ -- contradiction!
\end{proof}
\subsubsection{Analytic Luna slices}
The Luna Slice Theorem is an important theorem in invariant theory
of algebraic reductive groups (see \cite{Lun,Dre}). Here we
formulate a version of the Luna Slice Theorem for points over
local fields. For Archimedean F this was done by Luna himself in
[Lun75].

\begin{definition}
Let a reductive group $G$ act on an affine variety $X$. Let $\pi:
X(F) \to (X/G)(F)$ be the natural map. An open subset $U \subset
X(F)$ is called \textbf{saturated} if there exists an open subset
$V \subset (X/G)(F)$ such that $U = \pi^{-1}(V)$.
\end{definition}

\noindent We will use the following corollary of the Luna Slice
Theorem:

\begin{theorem}[\cite{AG_Gen_HC_RJR}, Theorem 2.3.17] \label{LocLuna}
Let a reductive group $G$ act on a smooth affine variety $X$. Let
$x \in X(F)$ be $G$-semisimple. Consider the natural action of the
stabilizer $G_x$ on the normal space $N_{Gx,x}^{X}$.
Then there exist\\
(i) an open $G(F)$-invariant $B$-analytic neighborhood $U$ of
$G(F)x$ in $X(F)$ with a
$G$-equivariant $B$-analytic retract $p:U \to G(F)x$ and\\
(ii) a $G_x$-equivariant $B$-analytic embedding $\psi:p^{-1}(x)
\hookrightarrow N_{Gx,x}^{X}(F)$ with an open saturated image such
that $\psi(x)=0$.
\end{theorem}

\begin{definition}
In the notation of the previous theorem, denote $S:= p^{-1}(x)$
and $N:=N_{Gx,x}^{X}(F)$. We call the quintuple $(U,p,\psi,S,N)$
an \textbf{analytic Luna slice at $x$}.
\end{definition}

\begin{corollary} \label{LocLunCor}
In the notation of the previous theorem, let $y\in p^{-1}(x)$.
Denote $z:=\psi(y)$. Then\\
(i) $(G(F)_x)_z=G(F)_y$\\
(ii) $N_{G(F)y,y}^{X(F)} \cong N_{G(F)_x z, z}^{N}$ as
$G(F)_y$-spaces\\
(iii) $y$ is $G$-semisimple if and only if $z$ is
$G_x$-semisimple.
\end{corollary}

\subsection{Generalized Harish-Chandra descent}
\label{GHC}
In this subsection we will prove the following theorem.
\begin{theorem} \label{Gen_HC}
Let a reductive group $G$ act on a smooth affine variety $X$. Let
$\chi$ be a character of $G(F)$. Suppose that for any
$G$-semisimple $x \in X(F)$ we have
$$\cD(N_{Gx,x}^X(F))^{G(F)_x,\chi}=0.$$ Then
$$\cD(X(F))^{G(F),\chi}=0.$$
\end{theorem}

\begin{remark}
In fact, the converse is also true. We will not prove it since we
will not use it.
\end{remark}

For the proof of this theorem we will need the following lemma

\begin{lemma}
Let a reductive group $G$ act on a smooth affine variety $X$. Let
$\chi$ be a character of $G(F)$. Let $U \subset X(F)$ be an open
saturated subset. Suppose that $\cD(X(F))^{G(F),\chi}=0.$ Then
$\cD(U)^{G(F),\chi}=0.$
\end{lemma}
\begin{proof}

Consider the quotient $X/G$. It is an affine algebraic variety.
Embed it in an affine space $\mathbb{A}^n$. This defines a map
$\pi:X(F) \to F^n$. Since $U$ is saturated, there exists an open
subset $V \subset (X/G)(F)$ such that $U=\pi^{-1}(V)$. Clearly
there exists an open subset $V' \subset F^n$ such that $V' \cap
(X/G)(F) = V$.

Let $\xi \in \cD(U)^{G(F),\chi}$. Suppose that $\xi$ is non-zero.
Let $x \in \Supp \xi$ and let $y:=\pi(x)$. Let $g \in
C^{\infty}_c(V')$ be such that $g(y) = 1$. Consider $\xi' \in
\cD(X(F))$ defined by $\xi'(f) := \xi(f \cdot (g \circ \pi))$.
Clearly, $\Supp(\xi') \subset U$ and hence we can interpret $\xi'$
as an element in $\cD(X(F))^{G(F),\chi}$. Therefore $\xi'=0$. On
the other hand, $x \in \Supp(\xi')$. Contradiction.
\end{proof}

\begin{proof}[Proof of Theorem \ref{Gen_HC}.] Let $x$
be a $G$-semisimple element. Let ($U_x$,$p_x$,$\psi_x$,
$S_x$,$N_x$) be an analytic Luna slice at $x$.

Let $\xi' = \xi |_{U_x}$. Then $\xi' \in \cD(U_x)^{G(F),\chi}.$ By
Frobenius descent it corresponds to $\xi'' \in
\cD(S_x)^{G_x(F),\chi}$.

The distribution $\xi''$ corresponds to a distribution $\xi''' \in
\cD(\psi_x(S_x))^{G_x(F),\chi}.$

However, by the previous lemma the assumption implies that
$\cD(\psi_x(S_x))^{G_x(F),\chi}=0.$ Hence $\xi'=0$.

Let $S \subset X(F)$ be the set of all $G$-semisimple points. Let
$U= \bigcup_{x\in S} U_x$. We saw that $\xi|_U=0$. On the other
hand, $U$ includes all the closed orbits, and hence by Corollary
\ref{OpenClosedAll} $U=X$.
\end{proof}

The following generalization of this theorem is proven in the same
way.

\begin{theorem} \label{Gen_HC_K}
Let a reductive group $G$ act on a smooth affine variety $X$. Let
$K \subset G(F)$ be an open subgroup and let $\chi$ be a character
of $K$. Suppose that for any $G$-semisimple $x \in X(F)$ we have
$$\cD(N_{Gx,x}^X(F))^{K_x,\chi}=0.$$ Then
$$\cD(X(F))^{K,\chi}=0.$$
\end{theorem}

\begin{remark}
If $K$ is an open B-analytic subgroup\footnote{In fact, any open
subgroup of a B-analytic group is B-analytic.} of $G(F)$ then the
theorem also holds for Schwartz distributions. Namely, if
$\Sc^*(N_{Gx,x}^X(F))^{K_x}=0$ for any $G$-semisimple $x \in X(F)$
then $\Sc^*(X(F))^{K}=0$. The proof is the same.
\end{remark}

The following generalization of Theorem \ref{Gen_HC} is proven in
the same way.
\begin{theorem} \label{HC_Thm}
Let a reductive group $G$ act on smooth affine varieties $X$ and
$Y$. Let $\chi$ be a character of $G(F)$. Suppose that for any
$G$-semisimple $x \in X(F)$ we have
$$\Sc^*((N_{Gx,x}^X\times Y)(F))^{G(F)_x,\chi}=0.$$ Then $\Sc^*(X(F)\times Y(F))^{G(F)_x,\chi}=0.$
\end{theorem}

\subsection{A stronger version} \label{sGHC}
In this subsection we provide means to validate the conditions of
Theorems \ref{Gen_HC} and \ref{Gen_HC_K} based on an inductive
argument.

More precisely, the goal of this subsection is to prove the
following theorem.

\begin{theorem} \label{Strong_HC}
Let a reductive group $G$ act on a smooth affine variety $X$. Let
$K \subset G(F)$ be an open subgroup and let $\chi$ be a character
of $K$. Suppose that for any $G$-semisimple $x \in X(F)$ such that
$$\cD(R_{G_x}(N_{Gx,x}^X))^{K_x,\chi}=0$$ we have
$$\cD(Q_{G_x}(N_{Gx,x}^X))^{K_x,\chi}=0.$$ Then for any $G$-semisimple $x \in X(F)$ we have $$\cD(N_{Gx,x}^X(F))^{K_x,\chi}=0.$$
\end{theorem}

 Together with Theorem \ref{Gen_HC_K}, this theorem gives the
following corollary.
\begin{corollary} \label{Strong_HC_Cor}
Let a reductive group $G$ act on a smooth affine variety $X$. Let
$K \subset G(F)$ be an open subgroup and let $\chi$ be a character
of $K$. Suppose that for any $G$-semisimple $x \in X(F)$ such that
$$\cD(R(N_{Gx,x}^X))^{K_x,\chi}=0$$ we have
$$\cD(Q(N_{Gx,x}^X))^{K_x,\chi}=0.$$ Then $\cD(X(F))^{K,\chi}=0.$
\end{corollary}

From now till the end of the subsection we fix $G$, $X$, $K$ and
$\chi$. Let us introduce several definitions and notation.

\begin{notation}
Denote \itemize  \item  $T \subset X(F)$ the set of all
$G$-semisimple points.
\item For $x,y \in T$ we say that $x>y$ if $G_x\supsetneqq G_y$.
\item $T_0 := \{x \in T\ |\ \cD(Q(N_{Gx,x}^X))^{K_x,\chi}=0 \} = \{x \in T\ |\ \cD((N_{Gx,x}^X))^{K_x,\chi}=0 \} .$
\end{notation}

\begin{proof}[Proof of Theorem \ref{Strong_HC}]
We have to show that $T=T_0$. Assume the contrary.

Note that every chain in $T$ with respect to our ordering has a
minimum. Hence by Zorn's lemma every non-empty set in $T$ has a
minimal element. Let $x$ be a minimal element of $T - T_0$. To get
a contradiction, it is enough to show that
$\cD(R(N_{Gx,x}^X))^{K_x,\chi}=0$.

Denote $R:=R(N_{Gx,x}^X)$. By Theorem \ref{Gen_HC_K}, it is enough
to show that for any $y\in R$ we have
$$\cD(N_{G(F)_x y,y}^R)^{(K_x)_y,\chi}=0.$$

Let $(U,p,\psi,S,N)$ be an analytic Luna slice at $x$.

Since $\psi(S)$ is open and contains 0, we can assume, upon
replacing $y$ by $\lambda y$ for some $\lambda \in F^{\times}$,
that $y \in \psi(S)$. Let $z \in S$ be such that $\psi(z)=y$. By
Corollary \ref{LocLunCor}, $G(F)_z=(G(F)_x)_y\subsetneqq G(F)_x$
and $N_{G(F)_x y,y}^R \cong N_{G z,z}^X(F)$. Hence $(K_x)_y = K_z$
and therefore
$$\cD(N_{G(F)_x y,y}^R)^{(K_x)_y,\chi} \cong \cD(N_{G
z,z}^X(F))^{K_z,\chi}.$$

However $z < x$ and hence $z \in T_0$ which means that $\cD(N_{G
z,z}^X(F))^{K_z,\chi}=0$.
\end{proof}

\begin{remark}
One can rewrite this proof such that it will use Zorn's lemma for
finite sets only, which does not depend on the axiom of choice.
\end{remark}

\begin{remark}
As before, Theorem \ref{Strong_HC} and Corollary
\ref{Strong_HC_Cor} also hold  for Schwartz distributions, with a
similar proof.
\end{remark}
\section{Localization Principle}\label{SecRedLocPrin}
In this section we formulate and prove localization principle. For
non-Archimedean local fields, a more general theorem is proven in
\cite[Section 1.5]{Ber}. The proofs we give in this section work
for local fields of characteristic zero. This section is based on
\cite[Appendix D]{AG_Gen_HC_RJR}.


\begin{theorem}[Localization Principle] \label{LocPrin2}
Let a reductive group $G$ act on a smooth algebraic variety $X$.
Let $Y$ be an algebraic variety and $\phi:X \to Y$ be an affine
algebraic $G$-invariant map. Let $\chi$ be a character of $G(F)$.
Suppose that for any $y \in Y(F)$ we have
$\cD_{X(F)}((\phi^{-1}(y))(F))^{G(F),\chi}=0$. Then
$\cD(X(F))^{G(F),\chi}=0$.
\end{theorem}

\begin{proof}
Clearly, it is enough to prove the theorem for the case when $X$
and $Y$ are affine. In this case there exists a categorical
quotient $X/G$ (see Theorem \ref{Quotient}). By the definition of
categorical quotient $\phi$ is a composition of $\pi_X:X \to X/G$
and some morphism $X/G \to Y$. Since every fiber of $\pi_X$ is
included in some fiber of $Y$,  it is enough to prove the theorem
for the case when $Y = X/G$ and $\phi = \pi_X$.

By the Generalized Harish-Chandra Descent (Corollary
\ref{Strong_HC_Cor}), it is enough to prove that for any
$G$-semisimple $x\in X(F)$, we have
$$\cD_{N_{Gx,x}^X(F)}(\Gamma(N_{Gx,x}^X))^{G_x(F),\chi} = 0.$$

Let $(U,p,\psi,S,N)$ be an analytic Luna slice at $x$. Clearly,
$$\cD_{N_{Gx,x}^X(F)}(\Gamma(N_{Gx,x}^X))^{G_x(F),\chi} \cong
\cD_{\psi(S)}(\Gamma(N_{Gx,x}^X))^{G_x(F),\chi}\cong
\cD_{S}(\psi^{-1}(\Gamma(N_{Gx,x}^X)))^{G_x(F),\chi}.$$ By
Frobenius descent (Theorem \ref{Frob}),
$$\cD_{S}(\psi^{-1}(\Gamma(N_{Gx,x}^X)))^{G_x(F),\chi}=\cD_{U}(G(F)\psi^{-1}(\Gamma(N_{Gx,x}^X)))^{G(F),\chi}.$$

By Lemma \ref{Gamma}, $$G(F)\psi^{-1}(\Gamma(N_{Gx,x}^X))= \{y \in
X(F)| x \in \overline{G(F)y} \}.$$
Hence by Corollary \ref{EquivClassClosed},
$G(F)\psi^{-1}(\Gamma(N_{Gx,x}^X))$ is closed in $X(F)$. Hence
$$\cD_{U}(G(F)\psi^{-1}(\Gamma(N_{Gx,x}^X)))^{G(F),\chi}
 =
\cD_{X(F)}(G(F)\psi^{-1}(\Gamma(N_{Gx,x}^X)))^{G(F),\chi}.$$ Now,
$$G(F)\psi^{-1}(\Gamma(N_{Gx,x}^X)) \subset
\pi_X(F)^{-1}(\pi_X(F)(x))$$ and we are given
$$\cD_{X(F)}(\pi_X(F)^{-1}(\pi_X(F)(x)))^{G(F),\chi}=0$$ for any
$G$-semisimple $x$.
\end{proof}

\begin{remark} \label{LocPrinS}
An analogous statement holds for Schwartz distributions and the
proof is the same.
\end{remark}

\begin{corollary} \label{LocPrinSub}
Let a reductive group $G$ act on a smooth algebraic variety $X$.
Let $Y$ be an algebraic variety and $\phi:X \to Y$ be an affine
algebraic $G$-invariant submersion. Suppose that for any $y \in
Y(F)$ we have $\Sc^*(\phi^{-1}(y))^{G(F),\chi}=0$. Then
$\Sc^*(X(F))^{G(F),\chi}=0$.
\end{corollary}

\begin{proof}
For any $y \in Y(F)$, denote $X(F)_y:=(\phi^{-1}(y))(F)$. Since
$\phi$ is a submersion, for any $y \in Y(F)$ the set $X(F)_y$ is a
smooth manifold. Moreover, $d\phi$ defines an isomorphism between
$N_{X(F)_y,z}^{X(F)}$ and $T_{Y(F),y}$ for any $z \in X(F)_y$.
Hence the bundle $CN_{X(F)_y}^{X(F)}$ is a trivial
$G(F)$-equivariant bundle.

We know that $$\Sc^*(X(F)_y)^{G(F),\chi}=0.$$ Therefore for any
$k$, we have
$$\Sc^*(X(F)_y,\Sym^k(CN_{X(F)_y}^{X(F)}))^{G(F),\chi}=0.$$ Thus by
Theorem \ref{NashFilt}, $$\Sc^*_{X(F)}(X(F)_y)^{G(F),\chi}=0.$$
Now, the Localization Principle implies that
$$\Sc^*(X(F))^{G(F),\chi}=0.$$
\end{proof}

\begin{remark} \label{RemLocVectSys}
Theorem \ref{LocPrin2} and Corollary \ref{LocPrinSub} admit
obvious generalizations to constant vector systems. The same
proofs hold.
\end{remark}

\section{Distributions versus Schwartz distributions}
This section is relevant only for Archimedean $F$.

In this section we show that if there are no $G(F)$-equivariant
Schwartz distributions on $X(F)$ then there are no
$G(F)$-equivariant distributions on $X(F)$. This section is based
on \cite{AG_Gen_HC_RJR}.

\begin{theorem}   \label{NoSNoDist}
Let a reductive group $G$ act on a smooth affine variety $X$. Let
$V$ be a finite-dimensional algebraic representation of $G(F)$.
Suppose that $$\Sc^*(X(F),V)^{G(F)}=0.$$ Then
$$\cD(X(F),V)^{G(F)}=0.$$
\end{theorem}

For the proof we will need the following definition and theorem.

\begin{definition}
(i) Let a topological group $K$ act on a topological space $M$. We
call a closed $K$-invariant subset $C \subset M$ \textbf{compact
modulo $K$} if there exists a compact subset $C' \subset M$ such
that $C \subset KC'.$

(ii) Let a Nash group $K$ act on a Nash manifold $M$. We call a
closed $K$-invariant subset $C \subset M$ \textbf{Nashly compact
modulo $K$} if there exist a compact subset $C' \subset M$ and
semi-algebraic closed subset $Z \subset M$ such that $C \subset Z
\subset KC'.$
\end{definition}

\begin{remark}
Let a reductive group $G$ act on a smooth affine variety $X$. Let
$K:=G(F)$ and $M:=X(F)$. Then it is easy to see that the notions
of compact modulo $K$ and Nashly compact modulo $K$ coincide.
\end{remark}

\begin{theorem} \label{CompSchwartz}
Let a Nash group $K$ act on a Nash manifold $M$. Let $E$ be a
$K$-equivariant Nash bundle over $M$. Let $\xi \in \cD(M,E)^K$ be
such that $\Supp(\xi)$ is Nashly compact modulo $K$. Then $\xi \in
\Sc^*(M,E)^K$.
\end{theorem}
For the proof see the next subsection.

\begin{proof}[Proof of Theorem  \ref{NoSNoDist}]
Fix any $y \in (X/G)(F)$ and denote $M:=\pi_X^{-1}(y)(F)$.

By the Localization Principle (Theorem \ref{LocPrin2} and Remark
\ref{RemLocVectSys}), it is enough to prove that
$$\Sc^*_{X(F)}(M,V)^{G(F)}=\cD_{X(F)}(M,V)^{G(F)}.$$ Choose $\xi
\in \cD_{X(F)}(M,V)^{G(F)}$. $M$ has a unique closed stable
$G$-orbit and hence a finite number of closed $G(F)$-orbits. By
Theorem \ref{CompSchwartz}, it is enough to show that $M$ is
Nashly compact modulo $G(F)$. Clearly $M$ is semi-algebraic.
Choose representatives $x_i$ of the closed $G(F)$-orbits in $M$.
Choose compact neighborhoods $C_i$ of $x_i$. Let $C' := \bigcup
C_i$. By Corollary \ref{OpenClosedAll}, $G(F)C' \supset M$.
\end{proof}

\subsection{$K$-invariant distributions compactly supported modulo
$K$.} \label{KInvAreSchwartz}

In this subsection we prove Theorem \ref{CompSchwartz}.

For the proof we will need the following lemmas.

\begin{lemma} \label{RelComp}
Let $M$ be a Nash manifold. Let $C \subset M$ be a compact subset.
Then there exists a relatively compact open (semi-algebraic)
subset $U \subset M$ that includes $C$.
\end{lemma}

\begin{proof}
For any point $x \in C$ choose an affine chart, and let $U_x$ be
an open ball with center at $x$ inside this chart. Those $U_x$
give an open cover of $C$. Choose a finite subcover
$\{U_i\}_{i=1}^n$ and let $U:= \bigcup_{i=1}^n U_i$.
\end{proof}

\begin{lemma} \label{takayata}
Let $M$ be a Nash manifold. Let $E$ be a Nash bundle over $M$. Let
$U \subset M$ be a relatively compact open (semi-algebraic)
subset. Let $\xi \in \cD(M,E)$. Then $\xi|_U \in \Sc^*(U,E|_U)$.
\end{lemma}

\begin{proof}
It follows from the fact that extension by zero $ext:\Sc(U,E|_U)
\to C_c^{\infty}(M,E)$ is a continuous map.
\end{proof}

\begin{proof}[Proof of Theorem \ref{CompSchwartz}]
Let $Z \subset M$ be a semi-algebraic closed subset and $C \subset
M$ be a compact subset such that $Supp(\xi) \subset Z \subset KC$.

Let $U \supset C$ be as in Lemma \ref{RelComp}. Let $\xi' :=
\xi|_{KU}$. Since $\xi|_{M - Z}=0$, it is enough to show that
$\xi'$ is Schwartz.

Consider the surjective submersion $m_U:K \times U \to KU$. Let
$$\xi'':=m_U^* (\xi') \in \cD(K \times U,m_U^*(E))^K.$$ By
Proposition \ref{EnoughPull}, it is enough to show that $$\xi''
\in \Sc^*(K \times U,m_U^*(E)).$$
By Frobenius descent, $\xi''$ corresponds to $\eta \in \cD(U,E)$.
It is enough to prove that $\eta \in \Sc^*(U,E)$. Consider the
submersion $m:K \times M \to M$ and let $$\xi''' := m^*(\xi) \in
\cD(K \times M,m^*(E)).$$ By Frobenius descent, $\xi'''$
corresponds to $\eta' \in \cD(M,E)$. Clearly $\eta = \eta'|_U$.
Hence by Lemma \ref{takayata}, $\eta \in \Sc^*(U,E)$.
\end{proof}

\section{Fourier transform and Homogeneity Theorems}\label{SecFour}
This section is based on \cite{AG_Gen_HC_RJR}.

Let $G$ be a reductive group   and $V$ be a finite-dimensional
 $F$-rational representation of $G$. Let $\chi$ be a
character of $G(F)$. In this section we provide some tools to
verify that $\Sc^*(Q(V))^{G(F),\chi}=0$ provided that
$\Sc^*(R(V))^{G(F),\chi}=0$.

\subsection{Formulation}
For this subsection let $B$ be a non-degenerate bilinear form on a
finite-dimensional vector space $V$ over $F$. We also fix an
additive character $\kappa$ of $F$. If $F$ is Archimedean we take
$\kappa(x):=e^{2\pi \mathrm{i} \re(x)}$.

\begin{notation}
We identify $V$ and $V^*$ via $B$ and endow $V$ with the self-dual
Haar measure with respect to $\psi$. Denote by $\Fou_B: \Sc^*(V)
\to \Sc^*(V)$ the Fourier transform. For any B-analytic manifold
$M$ over $F$ we also denote by $\Fou_B:\Sc^*(M \times V) \to
\Sc^*(M \times V)$ the partial Fourier transform.
\end{notation}

\begin{notation}
Consider the homothety action of $F^{\times}$ on $V$ given by
$\rho(\lambda)v:= \lambda^{-1}v$. It gives rise to an action
$\rho$ of $F^{\times}$ on $\Sc^*(V)$.

Let $|\cdot |$ denote the normalized absolute value. Recall that
for $F =\R$, $|\lambda|$ is equal to the classical absolute value
but for $F =\C$, $|\lambda| = (\re \lambda)^2+ (\Im \lambda)^2$.
\end{notation}

\begin{notation}
We denote by $\gamma(B)$ the Weil constant.  For its definition
see e.g. \cite[\S 2.3]{Gel}  for non-Archimedean $F$ and \cite[\S
1]{RS1} for Archimedean $F$.

For any $t\in F^{\times}$ denote $\delta_B(t)=
\gamma(B)/\gamma(tB)$.
\end{notation}

Note that $\gamma(B)$ is an 8-th root of unity and if $\dim V$ is
odd and $F \neq \C$ then $\delta_B$ is \textbf{not} a
multiplicative character.

\begin{notation}
We denote $$Z(B):=\{x \in V\ |\ B(x,x)=0 \}.$$
\end{notation}

\begin{theorem} [non-Archimedean homogeneity] \label{NonArchHom}
Suppose that $F$ is \textbf{non-Archimedean}. Let $M$ be a
B-analytic manifold over $F$. Let $\xi \in \Sc^*_{V\times
M}(Z(B)\times M)$ be such that $\Fou_B(\xi) \in \Sc^*_{V\times
M}(Z(B)\times M)$. Then for any $t \in F^{\times}$, we have
$\rho(t)\xi=\delta_B(t) |t|^{\dim V/2} \xi$ and $\xi=
\gamma(B)^{-1} \Fou_B(\xi)$. In particular, if $\dim V$ is odd
then $\xi = 0$.
\end{theorem}
For the proof see e.g. \cite[\S\S 8.1]{RS2} or \cite[\S\S
3.1]{JR}.

For the Archimedean version of this theorem we will need the
following definition.
\begin{definition}
Let $M$ be a B-analytic manifold over $F$. We say that a
distribution $\xi \in \Sc^*(V \times M)$ is
\textbf{adapted to $B$} if either \\
(i) for any $t \in F^{\times}$ we have $\rho(t)\xi=
\delta(t)|t|^{\dim V/2} \xi$ and $ \xi$ is proportional to $\Fou_B \xi$ or\\
(ii) $F$ is Archimedean and for any $t \in F^{\times}$ we have
$\rho(t)\xi= \delta(t)t|t|^{\dim V/2} \xi$.
\end{definition}

Note that if $\dim V$ is odd and $F \neq \C$ then every
$B$-adapted distribution is zero.

\begin{theorem} [Homogeneity Theorem] \label{ArchHom}
Let $M$ be a Nash manifold. Let $L\subset \Sc^*_{V\times
M}(Z(B)\times M)$ be a non-zero subspace such that for all $\xi
\in L $ we have $\Fou_B(\xi) \in L$ and $B \cdot \xi \in L$ (here
$B$ is viewed as a quadratic function).

Then there exists a non-zero distribution $\xi \in L$ which is
adapted to $B$.
\end{theorem}
For Archimedean $F$ we prove this theorem in subsection
\ref{AppRealHom}. For non-Archimedean $F$ it follows from Theorem
\ref{NonArchHom}.

We will also use the following trivial observation.

\begin{lemma}\label{FouInv}
Let a B-analytic group $K$ act linearly on $V$ and preserving $B$.
Let $M$ be a B-analytic $K$-manifold over $F$. Let $\xi \in
\Sc^*(V \times M)$ be a $K$-invariant distribution. Then
$\Fou_B(\xi)$ is also $K$-invariant.
\end{lemma}

\subsection{Applications}
The following two theorems easily follow form the results of the
previous subsection.

\begin{theorem} \label{Non_arch_Homog}
Suppose that $F$ is non-Archimedean. Let $G$ be a reductive group.
Let $V$ be a finite-dimensional  $F$-rational representation of
$G$. Let $\chi$ be character of $G(F)$. Suppose that
$\Sc^*(R(V))^{G(F),\chi}=0$. Let $V = V_1 \oplus V_2$ be a
$G$-invariant decomposition of $V$. Let $B$ be a $G$-invariant
symmetric non-degenerate bilinear form on $V_1$. Consider the
action $\rho$ of $F^{\times}$ on $V$ by homothety on $V_1$.

Then any $\xi \in \Sc^*(Q(V))^{G(F),\chi}$ satisfies
$\rho(t)\xi=\delta_B(t) |t|^{\dim V_1/2} \xi$ and $\xi= \gamma(B)
\Fou_B \xi$. In particular, if $\dim V_1$ is odd then $\xi = 0$.
\end{theorem}

\begin{theorem} \label{Homog}
Let $G$ be a reductive group. Let $V$ be a finite-dimensional
 $F$-rational representation of $G$. Let $\chi$ be character
of $G(F)$. Suppose that $\Sc^*(R(V))^{G(F),\chi}=0$. Let $Q(V) = W
\oplus (\bigoplus_{i=1}^kV_i)$ be a $G$-invariant decomposition of
$Q(V)$. Let $B_i$ be $G$-invariant symmetric non-degenerate
bilinear forms on $V_i$. Suppose that any $\xi \in
\Sc^*_{Q(V)}(\Gamma(V))^{G(F),\chi}$ which is adapted to each
$B_i$ is zero.

Then $\Sc^*(Q(V))^{G(F),\chi}=0.$
\end{theorem}
\begin{remark}
One can easily generalize Theorems \ref{Homog} and
\ref{Non_arch_Homog} to the case of constant vector systems.
\end{remark}

\subsection{Proof of the Archimedean Homogeneity Theorem}
\label{AppRealHom}
The goal of this subsection is to prove Theorem \ref{ArchHom} for
Archimedean $F$. We fix $V$ and $B$.

We will need some facts about the Weil representation. For a
survey on the Weil representation in the Archimedean case we refer
the reader to \cite[\S 1]{RS1}.

\begin{enumerate}
\item There exists a unique (infinitesimal) action
$\pi$ of $\sll_2(F)$ on $\Sc^*(V)$ such that\\
(i) $\pi(\begin{pmatrix}
  0 & 1 \\
  0 & 0
\end{pmatrix}) \xi = -\mathrm{i} \pi Re(B)
\xi$ and $\pi(\begin{pmatrix}
  0 & 0 \\
  -1 & 0
\end{pmatrix}) \xi = -\Fou_B^{-1} (\mathrm{i} \pi Re(B)
\Fou_B(\xi))$.\\
(ii) If $F=\C$ then $\pi(\begin{pmatrix}
  0 & \mathrm{i} \\
  0 & 0
\end{pmatrix})  = \pi(\begin{pmatrix}
  0 & 0 \\
  -\mathrm{i} & 0
\end{pmatrix})=0$

\item It can be lifted to an action of the metaplectic group
$\Mp(2,F)$.

We will denote this action by $\Pi$.

\item In case $F=\C$ we have $\Mp(2,F)=\SL_2(F)$ and in case $F=\R$ the
group $\Mp(2,F)$ is a connected 2-fold covering of $\SL_2(F)$. We
will denote by $\eps \in \Mp(2,F)$ the central element of order 2
satisfying $\SL_2(F)=\Mp(2,F)/\{1, \eps\}.$

\item In case $F=\R$ we have $\Pi(\eps)=(-1)^{\dim V}$ and therefore if $\dim
V$ is even then $\Pi$ factors through $\SL_2(F)$ and if $\dim V$
is odd then no nontrivial subrepresentation of $\Pi$ factors
through $\SL_2(F)$. In particular if $\dim V$ is odd then $\Pi$
has no nontrivial finite-dimensional representations, since every
finite-dimensional representation of $\Mp(2,F)$ factors through
$\SL_2(F)$. \label{FinDimSubrep}

\item In case $F=\C$ or in case $\dim V$ is even we have $\Pi(\begin{pmatrix}
  t & 0 \\
  0 & t^{-1}
\end{pmatrix}) \xi=\delta^{-1}(t)|t|^{-\dim V/2} \rho(t) \xi$ and $\Pi(\begin{pmatrix}
  0 & 1 \\
  -1 & 0
\end{pmatrix}) \xi=
\gamma(B)^{-1} \Fou_B \xi.$
\end{enumerate}

We also need the following straightforward lemma.

\begin{lemma} \label{sl2rep}
Let $(\Lambda, L)$ be a continuous finite-dimensional
representation of $\SL_2(\R)$.  Then there exists a non-zero $\xi
\in L$ such that either $$\Lambda(\begin{pmatrix}
  t & 0 \\
  0 & t^{-1}
\end{pmatrix})\xi=\xi \text{ and } \Lambda(\begin{pmatrix}
  0 & 1 \\
  -1 & 0
\end{pmatrix})\xi \text{ is proportional to }\xi$$ or
$$\Lambda(\begin{pmatrix}
  t & 0 \\
  0 & t^{-1}
\end{pmatrix})\xi=t \xi,$$
for all $t$.
\end{lemma}

Now we are ready to prove the theorem.
\begin{proof}[Proof of Theorem \ref{ArchHom}]
Without loss of generality assume $M=pt$.

Let $\xi \in L$ be a non-zero distribution. Let $L':=
U_{\C}(\sll_2(\R)) \xi \subset L$. Here, $U_{\C}$ means the
complexified universal enveloping algebra.

We are given that $\xi, \Fou_B(\xi) \in \Sc^*_V(Z(B))$. By Lemma
\ref{FinDim} below this implies that $L' \subset \Sc^*(V)$ is
finite-dimensional. Clearly, $L'$ is also a subrepresentation of
$\Pi$. Therefore by Fact (\ref{FinDimSubrep}), $F=\C$ or $\dim V$
is even. Hence $\Pi$ factors through $\SL_2(F)$.

Now by Lemma \ref{sl2rep} there exists $\xi' \in L'$ which is
$B$-adapted.
\end{proof}

\begin{lemma}\label{FinDim}
Let $V$ be a representation of $\sll_2$. Let $v \in V$ be a vector
such that $e^k v=f^n v=0$ for some $n,k$. Then the representation
generated by $v$ is finite-dimensional.\footnote{For our purposes
it is enough to prove this lemma for k=1.}
\end{lemma}
This lemma is probably well-known. Since we have not found any
reference we include the proof.
\begin{proof}
The proof is by induction on k.\\

\noindent Base k=1:\\
It is easy to see that $$e^l f^l v=l!(\prod_{i=0}^{l-1}(h-i))v$$
for all $l$. This can be checked by direct computation, and also
follows from the fact that $e^l f^l$ is of weight $0$, hence it
acts on the singular vector $v$ by its Harish-Chandra projection
which is $$\HC(e^l f^l)=l! \prod_{i=0}^{l-1}(h-i).$$

Therefore $(\prod_{i=0}^{n-1}(h-i))v=0$.

Hence $W:=U_{\C}(h) v$ is finite-dimensional and $h$ acts on it
semi-simply. Here, $U_{\C}(h)$ denotes the universal enveloping
algebra of $h$.  Let $\{v_i\}_{i=1}^m$ be an eigenbasis of $h$ in
$W$. It is enough to show that $U_{\C}(\sll_2)v_i$ is
finite-dimensional for any $i$. Note that $e|_W=f^n|_W=0$. Now,
$U_{\C}(\sll_2)v_i$ is finite-dimensional by the
Poincare-Birkhoff-Witt
Theorem.\\

\noindent Induction step:\\
Let $w:=e^{k-1}v$. Let us show that $f^{n+k-1} w=0$. Consider the
element $f^{n+k-1}e^{k-1}\in U_{\C}(\sll_2)$. It is of weight
$-2n$, hence by the Poincare-Birkhoff-Witt Theorem it can be
rewritten as a combination of elements of the form $e^a h^b f^c$
such that $c-a= n$ and hence $c \geq n$. Therefore
$f^{n+k-1}e^{k-1} v=0$.

Now let $V_1:=U_{\C}(\sll_2) v$ and $V_2:=U_{\C}(\sll_2) w$. By
the base of the induction $V_2$ is finite-dimensional, by the
induction hypotheses  $V_1/V_2$ is finite-dimensional, hence $V_1$
is finite-dimensional.
\end{proof}

\section{Singular support} \label{SingSupp}
In this section we assume that $F$ has zero characteristic. This
section is based on \cite{AG_AMOT}, \cite{Aiz} and \cite{Hef}.
%
%
In this section we will introduce the notion of singular support
of a distribution and prove some of its properties.

The most important property states that the singular support of a
distribution on $X$ is a (weakly) coisotropic subvariety of
$T^*X$. Therefore we dedicate the first subsection to coisotropic
and weakly coisotropic varieties.

\subsection{Coisotropic varieties} \label{coisotropic}

\begin{definition}
Let $X$ be a smooth algebraic variety and $\omega$ be a symplectic
form on it.
Let $Z\subset X$ be an algebraic subvariety. We call it {\bf $X$-coisotropic} if one of the following equivalent conditions holds.\\
(i) The ideal sheaf of regular functions that vanish on $\overline{Z}$ is closed under Poisson bracket. \\
(ii) At every smooth point $z \in Z$ we have  $T_zZ \supset
(T_zZ)^{\bot}$. Here, $(T_zZ)^{\bot}$ denotes the orthogonal
complement
to $T_zZ$ in $T_zX$ with respect to $\omega$. \\
(iii) For a generic smooth point $z \in Z$ we have $T_zZ \supset
(T_zZ)^{\bot}$.

If there is no ambiguity, we will call $Z$ a coisotropic variety.
\end{definition}
Note that every non-empty $X$-coisotropic variety is of dimension
at least $\frac{1}{2}\dim X$.

\begin{notation}
For a smooth algebraic variety $X$ we always consider the standard
symplectic form on $T^*X$. Also, we denote by $p_X:T^*X \to X$ the
standard projection.
\end{notation}

\begin{definition}
Let $(V,\omega)$ be a symplectic vector space with a fixed Lagrangian subspace $L \subset V$.
Let $p: V \to V/L$ be the standard projection. Let $W \subset V$ be a linear subspace.
We call it $V$-weakly coisotropic with respect to $L$ if one of the following equivalent conditions holds.\\
(i)$Z^ \bot \cap L \subset  Z\cap L  $\\
Here, $Z^{\bot}$ denotes the orthogonal complement with respect to
$\omega.$\\
(ii) $p(Z^ \bot) \subset p(Z)$.\\
(iii) $p(Z)^\bot \subset Z \cap L$.\\
Here, $p(Z)^\bot$ denotes the orthogonal complement in $L$ under
the identification $L \cong (V/L)^*$.
\end{definition}

\begin{definition}
Let $X$ be a smooth algebraic variety. Let $Z\subset T^*X$ be an algebraic subvariety. We call it {\bf $T^*X$-weakly coisotropic} if one of the following equivalent conditions holds.\\
(i)  At every smooth point $z \in Z$ the space $T_z(Z)$ is $T_z(T^*(X))$ -weakly coisotropic with respect to $Ker(d p_X).$\\
(ii)For a generic smooth point $z \in Z$ the space $T_z(Z)$ is $T_z(T^*(X))$ -weakly coisotropic with respect to $Ker(d p_X).$\\
(iii) For a generic smooth point  $x \in Z$ and for a generic
smooth point  $y \in p_X^{-1}(x) \cap Z$ we have
$$CN_{p_X(Z),x}^X \subset T_y(p_X^{-1}(x) \cap Z).$$
(iv) For any smooth point  $x \in p_X(Z)$ the fiber $p_X^{-1}(x)
\cap Z$ is locally invariant with respect to shifts by
$CN_{p_X(Z),x}^X$. I.e. for any point $y \in p_X^{-1}(x)$ the
intersection $$(y+CN_{p_X(Z),x}^X) \cap (p_X^{-1}(x)\cap Z)$$ is
Zariski open in $y+CN_{p_X(Z)}$.

If there is no ambiguity, we will call $Z$ a weakly coisotropic
variety.
\end{definition}
Note that every non-empty $T^*X$-weakly coisotropic variety is of
dimension at least $\dim X$.

The following lemma is straightforward.
\begin{lemma}
Any $T^*X$-coisotropic variety is $T^*X$-weakly coisotropic.
\end{lemma}
\begin{proposition}
Let $X$ be a smooth algebraic variety with a symplectic form on
it. Let $R \subset T^*X$ be an algebraic subvariety. Then there
exists a maximal $T^*X$-weakly coisotropic subvariety of $R$ i.e.
a $T^*X$-weakly coisotropic subvariety $T \subset R$ that includes
all $T^*X$-weakly coisotropic subvarieties of $R$.
\end{proposition}
\begin{proof}
Let $T'$ be the union of all smooth $T^*X$-weakly coisotropic
subvarieties of $R$. Let $T$ be the Zariski closure of $T'$ in
$R$. It is easy to see that $T$ is the maximal $T^*X$-weakly
coisotropic subvariety of $R$.
\end{proof}
The following lemma is trivial.
\begin{lemma}
Let $X$ be a smooth algebraic variety. Let a group $G$ act on $X$.
This induces an action on $T^*X$. Let $S \subset T^* X$ be a
$G$-invariant subvariety. Then the maximal $T^*X$-weakly
coisotropic subvariety of $S$ is also $G$-invariant.
\end{lemma}
\begin{notation}
Let $Y$ be a smooth algebraic variety. Let $Z \subset Y$ be a
smooth subvariety and $R \subset T^*Y$  be any subvariety. We
define {\bf the restriction $R|_Z \subset T^*Z$ of $R$ to $Z$} by
$R|_Z:=i^*(R)$, where $i:Z \to Y$ is the embedding.
\end{notation}
\begin{lemma} \label{Restriction}
Let $Y$ be a smooth algebraic variety. Let $Z \subset Y$ be a
smooth subvariety and $R \subset T^*Y$  be a weakly coisotropic
subvariety. Assume that any smooth point $z \in Z \cap p_Y(R)$ is
also a smooth point of $p_Y(R)$ and we have $T_z(Z \cap p_Y(R)) =
T_z(Z) \cap T_z (p_Y(R))$.

Then $R|_Z$ is $T^*Z$-weakly coisotropic.
\end{lemma}
\begin{proof}
Let $x \in Z$, let $M:= p_Y^{-1}(x) \cap R \subset p_Y^{-1}(x)$
and $L:=CN_{p_Y(R),x}^Y \subset p_Y^{-1}(x).$ We know that $M$ is
locally invariant with respect to shifts in $L$. Let $M':=
p_Z^{-1}(x) \cap R|_Z \subset p_Z^{-1}(x)$ and
$L':=CN_{p_Z(R|_Z),x}^Y \subset p_Z^{-1}(x).$ We want to show that
$M'$ is locally invariant with respect to shifts in $L'$. Let $q:
p_Y^{-1}(x) \to p_Z^{-1}(x)$ be the standard projection. Note that
$M'=q(M)$ and $L'=q(L)$. Now clearly $M'$ is locally invariant
with respect to shifts in $L'$.
\end{proof}

\begin{corollary} \label{PreGeoFrob}
Let $Y$ be a smooth algebraic variety.  Let an algebraic group $H$
act on $Y$. Let $q:Y \to B$ be an $H$-equivariant morphism. Let $O
\subset B$ be an orbit. Consider the natural action of $G$ on
$T^*Y$ and let $R \subset T^*Y$ be an $H$-invariant subvariety.
Suppose that $p_Y(R) \subset q^{-1}(O)$. Let $x \in O$. Denote
$Y_x:= q^{-1}(x)$. Then

\itemize{ \item if $R$ is $T^*Y$-weakly coisotropic then $R|_{Y_x}$ is $T^*(Y_x)$-weakly coisotropic.}
\end{corollary}
\begin{corollary} \label{GeoFrob}
In the notation of the previous corollary, if $R|_{Y_x}$ has no
(non-empty) $T^*(Y_x)$-weakly coisotropic subvarieties then $R$
has no (non-empty) $T^*(Y)$-weakly coisotropic subvarieties.
\end{corollary}

\begin{remark}
The  results  on weakly coistropic varieties  that we presented
here have versions for coistropic varieties, see \cite[section
5.1]{AG_AMOT}.
\end{remark}

\subsection{Definition and properties of singular support}
\label{WFSS}


\begin{notation}
\item Let $M,N$ be (smooth, algebraic, etc) manifolds. Let $S \subset(T^*(N))$. Let $\phi: M \to N$ be a morphism.
We denote $\phi ^*(S):=d(\phi)^*(S \times_N M)$.
\item Let $V$ be a linear space. For a point $x=(v,\phi) \in V \times V^*$ we denote $\widehat{x}=(\phi,-v) \in V^* \times V$, similarly for subset $X \subset V \times V^*$ we define $\widehat{X}$. for a (smooth, algebraic, etc) manifold  and a subset $X  \subset  T^*(M \times V)$ we denote $\widehat{X}_V  \subset  T^*(M \times V^*)$ in a similar
way.
\item Let $B$ be a non-degenerate bilinear form on $V$.
This gives an identification between $V$ and $V^*$ and therefore,
by the previous notation, defines maps $F_B: V \times V \to V
\times V$ and $F_B: T^*M \times V \times V \to T^*M \times V
\times V$. If there is bo ambiguity we will denote it by $F_V$.
\end{notation}

\begin{definition}
Let $F$ be non-Archimedean. Let $U \subset F^n$ be an open subset
and $\xi \in \Sc^*(U)$ be a distribution. We say that $\xi$ is
smooth at $(x_0, v_0) \in T^*U$ if there are open neighborhoods
$A$ of $x_0$ and $B$ of $v_0$ such that for any $\phi \in \Sc(A)$
there is an $N_\phi > 0$ for which for any $\lambda \in F$
satisfying $\lambda
> N_\phi$ we have $\widehat{(\phi \xi)}|_{\lambda B} = 0$. The
complement in $T^*U$ of the set of smooth pairs $(x_0, v_0)$ of
$\xi$ is called the wave front set of $\xi$ and denoted by
$WF(\xi)$.
\end{definition}

\begin{remark}
Sometimes in the literature (cf. \cite{Hef}) the wave front set is
defined to be a subset of $T^*U -U \times 0$. In our notation this
subset will be $WF(\xi) -U \times 0$.
\end{remark}


We will use the following proposition.

%

\begin{proposition}[\cite{Aiz}, Proposition 2.1.19] \label{Ramiiso}
Let $F$ be non-Archimedean. Let $V, U \subset F^n$ be open subsets
and $f: V \to U$ be an analytic isomorphism. Then for any $\xi \in
\Sc^*(V)$ we have $WF(f^*(\xi)) = f^*(WF(\xi))$.
\end{proposition}

\begin{corollary}
Let $F$ be non-Archimedean. Let $X$ be an $F$-analytic manifold.
Then we can define the the wave front set of any element in
$\Sc^*(X)$.
\end{corollary}

\begin{definition}
Let $X$ be a smooth algebraic variety let $\xi \in \Sc^*(X(F))$.
We will now define the singular support of $\xi$, it is an
algebraic subvariety of $T^*X$ and we will denote it by $SS(\xi)$.

In the case when $F$ is non-Archimedean we define it to be the
Zariski closure of $WF(\xi)$. In the case when $F$ is Archimedean
we define it to be the singular support of the $D_X$-module
generated by $\xi$.
\end{definition}

\begin{remark}
For Archimedean $F$ one can also define the wave front set in a
way similar to the above. In this case the singular support always
includes the Zariski closure of $WF(\xi)$ but is sometimes bigger.
\end{remark}
For readers who are not familiar with the theory of D-modules, the
next theorem summarizes the only properties of singular support
that we use. We will give more details on the algebraic theory of
D-modules in subsection \ref{AppDmod}. For a good introduction to
the algebraic theory of D-modules we refer the reader to
\cite{BerDmod} and \cite{Bor}.

\begin{theorem} \label{ThmSingSupp}
Let $X$ be a smooth algebraic variety. Then
\begin{enumerate}
\item \label{Supp2SS}
Let $\xi \in \Sc^*(X (F))$.  Then $\overline{\Supp(\xi)}_{Zar} =
p_X(SS(\xi))(F)$, where $\overline{\Supp(\xi)}_{Zar} $ denotes the
Zariski  closure of
 $\Supp(\xi)$.
\item \label{Ginv}
Let an algebraic group $G$ act on $X$. Let $\g$ denote the Lie
algebra of $G$. Let $\xi \in \Sc^*(X(F))^{G(F)}$. Then
$$SS(\xi) \subset \{(x,\phi) \in T^*X \, | \, \forall \alpha \in
\g \, \phi(\alpha(x)) =0\}.$$
\item \label{Fou}
Let $V$ be a linear space. Let $Z \subset X \times V$ be a closed
subvariety, invariant with respect to homotheties in $V$. Suppose
that $\Supp(\xi) \subset Z(F)$. Then $SS(\Fou_V(\xi)) \subset
F_V(p_{X \times V}^{-1}(Z))$.
\item \label{Gab}
Let $X$ be a smooth algebraic variety. Let $\xi \in \Sc^*(X(F))$.
Then $SS(\xi)$ is weakly coisotropic.
\end{enumerate}
\end{theorem}

For non-Archimedean fields these properties are proven in
\cite[Section 4]{Aiz}. For non-Archimedean fields, properties
\ref{Supp2SS} - \ref{Fou} are proven in subsection \ref{AppDmod}.
Property \ref{Gab}, which is crucial for us, follows from the
following theorem.

\begin{theorem}[Integrability Theorem] \label{Gaber}
Let $X$ be a smooth algebraic variety. Let ${\mathcal M}$ be a
finitely generated $D_X$-module. Then $SS({\mathcal M})$ is a
$T^*X$-coisotropic variety.
\end{theorem}
This is a special case of Theorem I in  \cite{Gab}. For similar
versions see also \cite{KKS, Mal}.

\begin{remark} Evidently, in the Archimedean case
the Integrability Theorem implies a stronger property. Namely it
implies that $SS(\xi)$ is coisotropic and not only weakly
coisotropic. We conjecture that this holds for the non-Archimedean
case also.
\end{remark}

\subsection{Distributions on non distinguished nilpotent orbits} \label{nondist}

\begin{definition}
Let $V$ be an algebraic finite dimensional representation over $F$
of a reductive group $G$. Suppose that there is a finite number of
$G$ orbits in $\Gamma(V).$ Let $x \in \Gamma(V).$ We will call it
$G$-distinguished, if $CN_{G x, x}^{Q(V)} \subset \Gamma(V^*)$. We
will call a $G$ orbit $G$-distinguished if all (or equivalently
one of) its elements are $G$- distinguished.

If there is no ambiguity we will omit the "$G$-".
\end{definition}

\begin{example} \label{group_case}
For the case of a semisimple group acting on its Lie algebra, the
notion of $G$-distinguished element coincides with the standard
notion of distinguished nilpotent element. In particular, in the
case when $G=SL_n$ and $V = sl_n$ the set of $G$-distinguished
elements is exactly the set of regular nilpotent elements.
\end{example}

\begin{proposition} \label{non_disting_no_cois}
Let $V$ be an algebraic finite dimensional representation over $F$
of a reductive group $G$. Suppose that there is a finite number of
$G$ orbits on $\Gamma(V).$ Let $W:=Q(V)$, let $A$ be the set of
non-distinguished elements in $\Gamma(V).$ Then there are no
non-empty $W \times W^*$-weakly coisotropic subvarieties of $A
\times \Gamma(V^*).$
\end{proposition}
The proof is clear.
\begin{corollary} \label{non_disting_no_dists}
Let $\xi \in \Sc^*(W)$ and suppose that  $\Supp(\xi) \subset
\Gamma(V)$ and $supp(\widehat{\xi}) \subset  \Gamma(V^*)$. Then
the set of distinguished elements in $\Supp(\xi)$ is dense in
$\Supp(\xi)$.
\end{corollary}
\begin{remark}
In the same way one can prove an analogous result for
distributions on $W \times M$ for any $B$-analytic manifold $M$.
\end{remark}

\subsection{D-modules and proof of properties \ref{Supp2SS} -
\ref{Fou} in the Archimedean case} \label{AppDmod}
%
In this subsection $X$ denotes a smooth affine variety defined
over $\R$. All the statements of this section extend automatically
to general smooth algebraic varieties defined over $\R$. In this
thesis we use only the case when $X$ is an affine space.

\begin{definition}
Let $D(X)$ denote the algebra of polynomial differential operators
on $X$.
We consider the filtration $F^{\leq i}D(X)$ on $D(X)$ given by the
order of differential operator.
\end{definition}

\begin{definition}
We denote by $\Gr D(X)$ the associated graded algebra of $D(X)$.

Define the symbol map $\sigma :D(X) \to \Gr D(X)$ in the following
way. Let $d \in D(X)$. Let $i$ be the minimal index such that $d
\in F^{\leq i}$. We define $\sigma(d)$ to be the image of $d$ in
 $(F^{\leq i}D(X)) / (F^{\leq i-1}D(X))$
\end{definition}

\begin{proposition}
$\Gr D(X) \cong {\mathcal O}(T^*X)$.
\end{proposition}
For proof see e.g. \cite{Bor}. 

\begin{notation}
Let $(V,B)$ be a quadratic space. \\
(i) We define a morphism of algebras $\Phi^D_V: D(X \times V) \to
D(X \times V)$ in the following way.

Consider $B$ as a map $B:V \to V^*$. For any $f \in V^*$ we set
$\Phi^D_V(f) :=
\partial_{B^{-1}(f)}$. For any $v \in V$ we set $\Phi^D_V(\partial _v) :=
-B(v)$ and for any $d \in D(X)$ we set $\Phi^D_V(d):=d$. \\
(ii) It defines a morphism of algebras $\Phi^O_V:{\mathcal
O}(T^*X) \to {\mathcal O}(T^*X)$.
\end{notation}

The following lemma is straightforward.
\begin{lemma}
Let $f$ be a homogeneous polynomial. Consider it as a differential
operator. Then $\sigma( \Phi^D_V(f))=\Phi^O_V(\sigma(f))$.
\end{lemma}

The D-modules we use in the paper are right D-modules. The
difference between right and left D-modules is not essential (see
e.g. section VI.3 in \cite{Bor}). We will use the notion of good
filtration on a D-module, see e.g. section II.4 in \cite{Bor}. Let
us now remind the definition of singular support of a module and a
distribution.

\begin{notation}
Let $M$ be a $D(X)$-module. Let $\alpha \in M$ be an element. Then
we denote by $Ann_{D(X)}$ the annihilator of $\alpha$.
\end{notation}

\begin{definition}
Let $M$ be a $D(X)$-module. Choose a good filtration on $M$.
Consider $grM$ as a module over $\Gr D(X) \cong {\mathcal
O}(T^*X)$. We define
$$SS(M):= \Supp (\Gr M) \subset T^*X.$$
This does not depend on the choice of the good filtration on $M$
(see e.g. \cite{Bor}, section II.4).

For a distribution $\xi \in S^*(X(\R))$ we define $SS(\xi)$ to be
the singular support of the module of distributions generated by
$\xi$.
\end{definition}

The following proposition is trivial.

\begin{proposition}
Let $I < D(X)$ be a right ideal. Consider the induced filtrations
on $I$ and $D(X)/I$. Then $\Gr (D(X)/I) \cong \Gr (D(X)) / \Gr
(I)$.
\end{proposition}

\begin{corollary}
Let $\xi \in S^*(X)$. Then $SS(\xi)$ is the zero set of $\Gr
(Ann_{D(X)} \xi)$.
\end{corollary}

\begin{corollary}
Let $I< {\mathcal O}(T^*X) $ be the ideal generated by
$\{\sigma(d) \, | \, d \in Ann_{D(X)}(\xi) \}$. Then $SS(\xi)$ is
the zero set of $I$.
\end{corollary}

\begin{corollary}
Property \ref{Ginv} holds.
\end{corollary}

\begin{lemma}
Let $\xi \in S^*(X)$. Let $Z \subset X$ be a closed subvariety
such that $\Supp(\xi) \subset Z(\R)$. Let $f \in {\mathcal O}(X)$
be a polynomial that vanishes on $Z$. Then there exists $k \in \N
$ such that $f^k\xi=0$.
\end{lemma}
\begin{proof}$ $

Step 1. Proof for the case when $X$ is affine space and $f$ is a
coordinate function.\\
This follows from the proof of Corollary 5.5.4 in \cite{AG1}.

Step 2. Proof for the general case.\\
Embed $X$ into an affine space $A^N$ such that $f$ will be a
coordinate function and consider $\xi$ as distribution on $A^N$
supported in $X$. By Step 1, $f^k\xi=0$ for some $k$.
\end{proof}

\begin{corollary}
Property \ref{Supp2SS} holds.
\end{corollary}

\begin{proposition}
Property \ref{Fou} holds. Namely:

Let $(V,B)$ be a quadratic space. Let $Z \subset X \times V$ be a
closed subvariety, invariant with respect to homotheties in $V$.
Suppose that $\Supp(\xi) \subset Z(\R)$. Then $SS(\Fou_V(\xi))
\subset F_V(p_{X \times V}^{-1}(Z))$.
\end{proposition}
\begin{proof}
Let $f \in {\mathcal O}(X\times V)$ be homogeneous with respect to
homotheties in $V$. Suppose that $f$ vanishes on $Z$. Then
$\Phi^D_V(f^k) \in Ann_{D(X)} (\Fou_V(\xi))$. Therefore $\sigma
(\Phi^D_V(f^k))$ vanishes on $SS(\Fou_V(\xi))$. On the other hand,
$\sigma(\Phi^D_V(f^k)) =
\Phi^O_V(\sigma(f^k))=(\Phi^O_V(\sigma(f)))^k$. Hence
$SS(\Fou_V(\xi))$ is included in the zero set of
$\Phi^O_V(\sigma(f))$. Intersecting over all such $f$ we obtain
the required inclusion.
\end{proof}

\chapter{Gelfand Pairs and Gelfand-Kazhdan criterion} \label{ChGK}
This chapter is based on \cite{GK,AGS1,AG_Gen_HC_RJR,AG_AMOT}.
\section{Gelfand Pairs}

\begin{definition}
 Let $G$ be a reductive group. By an \textbf{admissible
 representation of} $G$ we mean an admissible smooth representation of
 $G(F)$ if $F$ is non-Archimedean (see \cite{BZ}) and admissible
 smooth \Fre representation of $G(F)$ if $F$ is Archimedean (see section \ref{SmoothFre}).
 \end{definition}

 We now introduce three a-priori distinct notions of Gelfand pair.

 \begin{definition}\label{GPs}
 Let $H \subset G$ be a pair of reductive groups.
 \begin{itemize}
 \item We say that $(G,H)$ satisfy {\bf GP1} if for any irreducible
 admissible representation $(\pi,E)$ of $G$
 we have
 $$\dim \Hom_{H(F)}(E,\cc) \leq 1.$$

 \item We say that $(G,H)$ satisfy {\bf GP2} if for any irreducible
 admissible representation $(\pi,E)$ of $G$
 we have
 $$\dim \Hom_{H(F)}(E,\cc) \cdot \dim \Hom_{H}(\widetilde{E},\cc)\leq
 1.$$

 \item We say that $(G,H)$ satisfy {\bf GP3} if for any irreducible
 {\bf unitary} representation $(\pi,\mathcal{H})$ of $G(F)$ on a
 Hilbert space $\mathcal{H}$ we have
 $$\dim \Hom_{H(F)}(\mathcal{H}^{\infty},\cc) \leq 1.$$
 We will call pairs that satisfy GP1 \textbf{Gelfand pairs} and pairs that satisfy GP3 \textbf{generalized Gelfand
 pairs}.
 \end{itemize}

 \end{definition}

 Property GP1 was established by Gelfand and Kazhdan in certain
 $p$-adic cases (see \cite{GK}). Property GP2 was introduced in
 \cite{Gross} in the $p$-adic setting. Property GP3 was studied
 extensively by various authors both in the real and $p$-adic settings (see e.g.
 \cite{vD-P}, \cite{vD}, \cite{Bos-vD}).

 We have the following straightforward proposition.

 \begin{proposition}
 $GP1 \Rightarrow GP2 \Rightarrow GP3.$
 \end{proposition}

 \begin{remark}
 It is not known whether some of these notions are equivalent.
 \end{remark}

 The most powerful tool for proving Gelfand property is the
 following theorem.

 \begin{theorem}[Gelfand-Kazhdan criterion]\label{DistCrit}
 Let $H \subset G$ be reductive groups and let $\sigma$ be an
 involutive anti-automorphism of $G$ and assume that $\sigma(H)=H$.
 Suppose $\sigma(\xi)=\xi$ for all bi $H(F)$-invariant Schwartz
 distributions $\xi$ on $G(F)$. Then $(G,H)$ satisfies GP2.
 \end{theorem}
 For non-Archimedean $F$ this theorem is proven in \cite{GK}. We
 adapt their proof to the Archimedean case in the next subsection.

\begin{remark}
This theorem is the main reason for introducing the property
$GP2$.
\end{remark}

\begin{remark}
In \cite{Thomas} it is proven that if $\sigma$ preserves for all
bi $H(F)$-invariant positive definite
 distributions on $G(F)$ then $(G,H)$ satisfies GP3.
\end{remark}

In some cases, GP2 is known to be equivalent to GP1. For example,
see Corollary \ref{GKCor} below.

 \begin{theorem}\label{GKRep}
 Let $G$ be a reductive group and let $\sigma$ be an
 anti-involution of $G$. Let $\theta$ be the
 involution of $G$ defined by $\theta(g):=\sigma(g^{-1})$. Let
 $(\pi, E)$ be an irreducible admissible representation of $G$.
 Suppose that either $\sigma$ preserves all conjugacy classes in
 $G(F)$ or $\charact F=0$ and $\sigma$ preserves all semisimple conjugacy classes in
 $G(F)$.

 Then $\widetilde{E} \cong E^{\theta}$, where $\widetilde{E}$
 denotes the smooth contragredient representation and $E^{\theta}$
 is $E$ twisted by $\theta$.
 \end{theorem}
 \begin{proof}
Since the representations $\widetilde{E}$ and $E^{\theta}$ are
irreducible, they are isomorphic if and only if their characters
are identical (see e.g. \cite[Theorem 8.1.5]{Wal1}).

If $F$ is non-archimedean and $\sigma$ preserves all conjugacy
classes in $G(F)$ then the identity of characters follows from
Theorem \ref{Constructibility}. If $\charact F=0$ and $\sigma$
preserves all semisimple conjugacy classes in $G(F)$ then this
follows from \cite[Corollary 7.6.3]{AG_Gen_HC_RJR}.
 \end{proof}

 \begin{corollary} \label{GKCor}
 Let $H \subset G$ be reductive groups and let $\sigma$ be an
 anti-involution of $G$ that preserves all conjugacy classes and such that
 $\sigma(H)=H$. Then $GP1$ is equivalent to $GP2$ for the pair
 $(G,H)$.
 \end{corollary}
Choosing $\sigma(g):=g^t$ we obtain
 \begin{corollary} \label{GKCorWeakPairs}
$GP1$ is equivalent to $GP2$ for the pairs
 $(\GL_{n+1},\GL_{n})$ and $(O(V \oplus F),O(V))$.
 \end{corollary}

\section[Smooth \Fre representations and proof of Theorem \ref{DistCrit}]{Smooth \Fre representations and proof of Theorem \ref{DistCrit} in the Archimedean
case} \label{SmoothFre}

The theory of representations in the context of \Fre spaces is
developed in \cite{Cas2} and \cite{Wal2}. We present here a
well-known slightly modified version of that theory. In this
section the field $F$ is Archimedean.

\begin{definition}
Let $V$ be a complete locally convex topological vector space. A
representation $(\pi,V,G)$ is a continuous map $G \times V \to V$.
A representation is called \textbf{\Fre} if there exists a
countable family of semi-norms $\rho_{i}$ on $V$ defining the
topology of $V$ and such that the action of $G$ is continuous with
respect to each $\rho_{i}$. We will say that $V$ is \textbf{smooth
\Fre} representation if, for any $X \in \g$ the differentiation
map $ v \mapsto \pi(X)v$ is a continuous linear map from $V$ to
$V$.
\end{definition}

An important class of examples of smooth \Fre representations is
obtained from continuous Hilbert representations $(\pi,H)$ by
considering the subspace of smooth vectors $H^{\infty}$ as a \Fre
space (see \cite[section 1.6]{Wal1} and \cite[section
11.5]{Wal2}).

We will consider mostly smooth \Fre representations.

\begin{remark}
In the language of \cite{Wal2} and \cite{Cas} the representations
above are called smooth \Fre representations of moderate growth.
\end{remark}

Recall that a smooth \Fre representation is called
\textit{admissible} if its underlying $(\g,K)$-module is
admissible (in particular finitely generated). In what follows
\textit{admissible representation} will always refer to admissible
smooth \Fre representation.

For an admissible (smooth \Fre) representation $(\pi,E)$ we denote
by $(\widetilde{\pi},\widetilde{E})$ the contragredient
representation.

We will require the following corollary of the globalization
theorem of Casselman and Wallach (see \cite{Wal2} , chapter 11).
\begin{theorem}\label{cortoCW}
Let $E$ be an admissible \Fre representation, then there exists a
continuous Hilbert space representation $(\pi,H)$ such that
$E=H^{\infty}$.
\end{theorem}
This theorem follows easily from the embedding theorem of
Casselman combined with Casselman-Wallach globalization theorem.

\Fre representations of $G$ can be lifted to representations of
$\Sc(G)$, the Schwartz space of $G$. 

For a \Fre representation $(\pi,E)$ of $G$, the algebra $\Sc(G)$
acts on $E$ through
\begin{equation}\label{operator}
\pi(\phi)=\int_{G} \phi(g) \pi(g)dg
\end{equation}
(see \cite{Wal1}, section 8.1.1).\\
The following lemma is straightforward:

\begin{lemma}\label{generalShalika}
Let $(\pi,E)$ be an admissible \Fre representation of $G$ and let
$\lambda \in E^*$. Then $\phi \to \pi(\phi)\lambda$ is a
continuous map $\Sc(G) \to \widetilde{E}$.
\end{lemma}

The following proposition follows from Schur's lemma for $(\g,K)$
modules (see \cite{Wal1} page 80) in light of Casselman-Wallach
theorem.

\begin{proposition}\label{SchurCW}
Let $G$ be a real reductive group. Let $W$ be a \Fre
representation of $G$ and let $E$ be an irreducible admissible
representation of $G$. Let $T_{1},T_{2}:W \hookrightarrow E$ be
two embeddings of $W$ into $E$. Then $T_{1}$ and $T_{2}$ are
proportional.
\end{proposition}

Now we are ready to prove Theorem \ref{DistCrit}.

\begin{proof}[Proof of Theorem \ref{DistCrit}.]
Let $(\pi,E)$ be an irreducible admissible \Fre representation. If
$E$ or $\widetilde{E}$ are not distinguished by $H$ we are done.
Thus we can assume that there exists a non-zero $\lambda:E \to
\cc$ which is $H$-invariant. Now let $\ell_{1},\ell_{2}$ be two
non-zero $H$-invariant functionals on $\widetilde{E}$. We wish to
show that they are proportional. For this we define two
distributions $D_{1},D_{2}$ as follows
$$D_{i}(\phi)=\ell_{i}(\pi(\phi)\lambda)$$ for $i=1,2$.
Here $\phi \in \Sc(G)$. Note that $D_{i}$ are also Schwartz
distributions. Both distributions are bi-$H$-invariant and hence,
by the assumption, both distributions are $\sigma$ invariant. Now
consider the bilinear forms on $\Sc(G)$ defined by
$$B_{i}(\phi_{1},\phi_{2})=D_{i}(\phi_{1}*\phi_{2}).$$
Since $E$ is irreducible, the right kernel of $B_{1}$ is equal to
the right kernel of $B_{2}$. We now use the fact that $D_{i}$ are
$\sigma$ invariant. Denote by $J_{i}$ the left kernels of $B_{i}$.
Then $J_{1}=J_{2}$ which we denote by $J$. Consider the \Fre
representation $W=\Sc(G)/J$ and define the maps $T_{i}: \Sc(G) \to
\widetilde{\widetilde{E}} \cong E$ by
$T_{i}(\phi)=\pi(\phi)\ell_{i}$. These are well defined by Lemma
\ref{generalShalika} and we use the same letters to denote the
induced maps $T_{i}: W \to E$. By Proposition \ref{SchurCW},
$T_{1}$ and $T_{2}$ are proportional and hence $\ell_{1}$ and
$\ell_{2}$ are proportional and the proof is complete.
\end{proof}

\section{Strong Gelfand Pairs}
Like in the previous section, there are several notions of strong
Gelfand pair.

\begin{definition}\label{SGPs}
 Let $H \subset G$ be a pair of reductive groups.
 \begin{itemize}
 \item We say that $(G,H)$ that is a strong
Gelfand pair or that $(G,H)$ satisfy {\bf SGP1} if for any
irreducible
 admissible representations $\pi$ of $G$ and $\tau$ of $H$
 we have
 $$\dim \Hom_{H(F)}(\pi,\tau) \leq 1.$$

 \item We say that $(G,H)$ satisfy {\bf SGP2} if for any
irreducible
 admissible representations $\pi$ of $G$ and $\tau$ of $H$
 we have
 $$\dim \Hom_{H(F)}(\pi,\tau) \cdot \dim \Hom_{H}(\widetilde{\pi},\widetilde{\tau})\leq
 1.$$

 \item We say that $(G,H)$ is a generalized strong Gelfand pair or that $(G,H)$ satisfy {\bf SGP3} if for any irreducible
 {\bf unitary} representations $\pi$ of $G$ and $\tau$ of $H$ we have
 $$\dim \Hom_{H(F)}(\pi,\tau) \leq 1.$$
 \end{itemize}

 \end{definition}

We will use the following two known theorems.
\begin{theorem}\label{IrrProd}
Let $G$ and $H$ be reductive groups. Let $\pi$ and $\tau$ be
irreducible admissible smooth representations of $G$ and $H$
respectively. Then $\pi \otimes \tau$ is irreducible
representation $G \times H$.
\end{theorem}

Is theorem is proven by reduction to an analogous statement on
finite-dimensional representations of algebras with 1. For the
statement see \cite[Ch. VIII, \S 7, Proposition 8]{Bou}.

For the non-Archimedean case this reduction is done in
\cite[Proposition 2.16]{BZ}. For the Archimedean case see e.g.
\cite[Appendix A]{AG_AMOT}.

\begin{theorem}\label{HomToProd}
Let  $H$ be a reductive group. Let $\pi$ and $\tau$ be smooth
(\Fre) representations of $H$. Suppose that $\tau$ is admissible.

Then $Hom(\pi,\tau)$ is canonically isomorphic to $Hom(\pi \otimes
\widetilde{\tau},\C).$
\end{theorem}
\begin{proof}[Sketch of proof.]
Consider the canonical embeddings
$$ Hom(\pi,\tau) \hookrightarrow Hom(\pi \otimes
\widetilde{\tau},\C) \hookrightarrow
Hom(\pi,\widetilde{\widetilde{\tau}}).$$ Since for admissible
representations $\widetilde{\widetilde{\tau}} \cong \tau$, those
embeddings are isomorphisms.
\end{proof}

\begin{corollary}\label{DeltaToStrong}
Let $G$ and $H$ be real reductive groups. The pair $(G, H)$ is a
strong Gelfand pair (respectively satisfies SGP2) if and only if
the pair $(G\times H, \Delta H)$ is a Gelfand pair (respectively
satisfies GP2).
\end{corollary}
Using Corollary \ref{GKCor} we obtain

 \begin{corollary} \label{SGKCor}
 Let $H \subset G$ be reductive groups and let $\sigma$ be an
 anti-involution of $G$ that preserves all conjugacy classes and such that
 $\sigma(H)=H$.

 Then $SGP1$ is equivalent to $SGP2$ for the pair
 $(G,H)$. In particular $SGP1$ is equivalent to $SGP2$ for the pair
 $(\GL_{n+1},\GL_{n})$.
 \end{corollary}

Corollary \ref{DeltaToStrong} also enables to formulate the
following analog of Gelfand-Kazhdan criterion for strong Gelfand
pairs.

\begin{theorem}\label{StrongDistCrit}
 Let $H \subset G$ be reductive groups and let $\sigma$ be an
 involutive anti-automorphism of $G$ and assume that $\sigma(H)=H$.
 Suppose $\sigma(\xi)=\xi$ for all $AdH(F)$-invariant Schwartz
 distributions $\xi$ on $G(F)$. Then $(G,H)$ satisfies SGP2.
\end{theorem}
\begin{proof}
Define $\sigma':G \times H \to G\times H$ by
$\sigma'(g,h):=(\sigma(g),\sigma(h))$. Let $\Delta H < G \times H$
be the diagonal. Consider the projection $G \times H \to H$. By
Frobenius descent (Theorem \ref{Frob}), the assumption implies
that any $\Delta H$-bi-invariant distribution on $G \times H$ is
invariant with respect to $\sigma'$.

Hence by Theorem \ref{DistCrit} $(G \times H, \Delta H)$ satisfies
GP2 and hence by the previous corollary $(G,H)$ satisfies SGP2.
\end{proof}

\begin{corollary}\label{StrongDist2StrongRep}
Theorem \ref{StrongDist} implies Theorem \ref{StrongRep}.
\end{corollary}

\chapter[Weak Gelfand property]{Gelfand property for the pairs $(\GL_{n+1}(F),\GL_n(F))$ and
$(\mathrm{O}_{n+1}(F),\mathrm{O}_n(F))$}\label{ChWeak} For fields
of characteristic zero the results of this chapter are much weaker
then the main results. I decided to include this chapter for two
reasons. First, the proofs in this chapter work over any local
field, while the analogs of the main results for fields of
positive characteristic are not known. Second, the proofs in this
chapter are shorter and clearly show how to use some of the tools.
The first section of this chapter is based on \cite{AGS1} and the
second on \cite{AGS2}.

\section{The pair $(\mathrm{GL}_{n+1}(F),\mathrm{GL}_n(F))$}
Consider the standard imbedding $\mathrm{GL}_n(F) \hookrightarrow
\mathrm{GL}_{n+1}(F)$. We consider the two-sided action of
$\mathrm{GL}_n(F) \times \mathrm{GL}_n(F)$ on
$\mathrm{GL}_{n+1}(F)$ defined by
$(g_1,g_2)h:=g_1 hg_2^{-1}$. The goal of this section is to prove
the following theorem.
\begin{theorem} \label{WeakGL}
Any $\GL_n(F) \times \GL_n(F)$ invariant distribution on
$\GL_{n+1}(F)$ is invariant with respect to transposition.
\end{theorem}
By Theorem \ref{DistCrit} and Corollary \ref{GKCorWeakPairs}  this
implies

\begin{theorem}
The pair $(\GL_{n+1}(F),\GL_n(F))$ is a Gelfand pair.
\end{theorem}

Since any character of $\mathrm{GL}_{n}(F)$ can be extended to
$\mathrm{GL}_{n+1}(F)$, we obtain

\begin{corollary}\label{WeakGLCor}
Let $(\pi,E)$ be an irreducible admissible representation of
$\mathrm{GL}_{n+1}(F)$ and let $\chi$ be a character of
$\mathrm{GL}_{n}(F)$.  Then
$$\dim Hom_{\mathrm{GL}_n(F)}(\pi,\chi) \leq 1.$$
\end{corollary}

For the proof of Theorem \ref{WeakGL} we will use the following
notation.
\begin{notation} \label{NotObjects}
Denote $H:=H_n:=\mathrm{GL}_n$. Denote
$$G:=G_n:=\{(h_1,h_2)\in \mathrm{GL}_n \times \mathrm{GL}_n|
\, det(h_1)=det(h_2)\}.$$ We consider $H$ to be diagonally
embedded to $G$.

Consider the action of the 2-element group $S_2$ on $G$ given by
the involution $(h_1,h_2) \mapsto ({h_2^{-1}}^t,{h_1^{-1}}^t)$. It
defines a semidirect product $\widetilde{G} := \widetilde{G}_n :=
G \rtimes S_2$. Denote also $\tH:= \tH_n:= H_n \rtimes S_2$.

Let $V=F^n$ and $X:=X_n := \mathrm{gl}_{n} \times V \times V^*$.

The group $\tG$ acts on $X$ by
$$(h_1,h_2)(A,v,\phi):= (h_1Ah_2^{-1},h_1v, {h_2^{-1}}^t \phi) \text{ and }$$
$$\sigma(A,v,\phi):=(A^t,\phi^t,v^t)$$
where $(h_1,h_2) \in G$ and $\sigma$ is the generator of $S_2$.
Note that $\tG$ acts separately on $\mathrm{gl}_n$ and on $V
\times V^*$. Define a character $\chi$ of $\widetilde{G}$ by
$\chi(g,s):= sign(s)$.
\end{notation}

We will show that the following theorem implies Theorem
\ref{WeakGL}.
\begin{theorem} \label{WeakGL-main}
$\Sc^*(X(F))^{\tG(F),\chi}=0$.
\end{theorem}

\subsection{Proof that Theorem \ref{WeakGL-main} implies Theorem
\ref{WeakGL}} \label{SecReduction}
We will divide this reduction to several propositions.\\
Consider the action of $\tG_n$ on $\mathrm{GL}_{n+1}$ and on
$\mathrm{gl}_{n+1}$, where $G_n$ acts by the two-sided action and
the generator of $S_2$ acts by transposition.

\begin{proposition} \label{WeakGL-Red1}
If  $\cD(\mathrm{GL}_{n+1}(F))^{\tG_n(F),\chi}=0$ then Theorem
\ref{WeakGL} holds.
\end{proposition}
The proof is straightforward.

\begin{proposition}
If  $\Sc^*(\mathrm{GL}_{n+1}(F))^{\tG_n(F),\chi}=0$ then
$\cD(\mathrm{GL}_{n+1(F)})^{\tG_n(F),\chi}=0$.
\end{proposition}
Follows from Theorem \ref{NoSNoDist}.
\begin{proposition} \label{WeakGL-Red2}
If  $\Sc^*(\mathrm{gl}_{n+1}(F))^{\tG_n(F),\chi}=0$ then
$\Sc^*(\mathrm{GL}_{n+1}(F))^{\tG_n(F),\chi}=0$.
\end{proposition}
\begin{proof}\footnote{This proposition is an adaption of a
statement in \cite{Ber}, section 2.2.} Let $\xi \in
\Sc^*(\mathrm{GL}_{n+1}(F))^{\tG_n(F),\chi}$. We have to prove
$\xi=0$. Assume the contrary. Take $p \in \mathrm{Supp}(\xi)$. Let
$t=\mathrm{det}(p)$. Let $f\in C_c^{\infty}(F)$ be such that $f$
vanishes in a neighborhood of zero and $f(t) \neq 0$. Consider the
determinant map $\mathrm{det}:\mathrm{GL}_{n+1}(F) \to F$.
Consider $\xi':=(f \circ \mathrm{det})\cdot \xi$. It is easy to
check that $\xi' \in \Sc^*(\mathrm{GL}_{n+1}(F))^{\tG_n(F),\chi}$
and $p \in \mathrm{Supp}(\xi')$. However, we can extend $\xi'$ by
zero to $\xi'' \in \Sc^*(\mathrm{gl}_{n+1}(F))^{\tG_n(F),\chi}$,
which is zero by the assumption. Hence $\xi'$ is also zero.
Contradiction.
\end{proof}
%
\begin{proposition} \label{WeakGL-Red4}
If $\Sc^*(X_{n}(F))^{\tG_n(F),\chi}=0$ then
$\Sc^*(\mathrm{gl}_{n+1}(F))^{\tG_n(F),\chi}=0$.
\end{proposition}
\begin{proof}
Note that $\mathrm{gl}_{n+1}$ is isomorphic as a
$\tG_n$-equivariant space to $X_n \times F$ where the action on
$F$ is trivial. This isomorphism is given by
$$ \left(
  \begin{array}{cc}
    A_{n \times n} & v_{n\times 1} \\
    \phi_{1\times n} & t\\
  \end{array}
\right) \mapsto ((A , v ,\phi),t) .$$ The proposition now follows
from Proposition \ref{Product}.
\end{proof}
This finishes the proof that Theorem \ref{WeakGL-main} implies
Theorem \ref{WeakGL}.
\subsection{Proof of Theorem \ref{WeakGL-main}}
We will now stratify $X (=gl_n \times V \times V^*)$ and deal with
each strata separately.
\begin{notation}
Denote $W:=W_n:=V_n \oplus V_n^*.$ Denote by $Q^i:=Q^i_n \subset
\mathrm{gl}_n$ the set of all matrices of rank $i$. Denote
$Z^i:=Z^i_n:=Q^i_n \times W_n$.
\end{notation}
Note that $X = \bigcup Z^i$. Hence by Theorems \ref{Filt_nonarch}
and \ref{NashFilt}, it is enough to prove the following
proposition.

\begin{proposition} \label{Enough}
$ $\\
(i) If $F$ is non-Archimedean then $\Sc^*(Z^i(F))^{\tG(F),\chi}=0$ for any $i$.\\
(ii) If $F$ is Archimedean then $$\Sc^*(Z^i(F), Sym^k(CN_{Z^i}^X
)(F))^{\tG(F),\chi}=0$$ for any $k$ and $i$.
\end{proposition}

We will use the following important lemma.
\begin{lemma} \label{keyGLWeak}
$\Sc^*(W(F))^{\tH(F),\chi}=0.$
\end{lemma}
For proof see subsection \ref{ProofKeyGLWeak} below.

\begin{corollary}
Proposition \ref{Enough} holds for $i=n$.
\end{corollary}
\begin{proof}
Clearly, one can extend the actions of $\widetilde{G}$ on $Q^n$
and on $Z^n$ to actions of $\widetilde{GL_n \times GL_n}:= (GL_n
\times GL_n) \rtimes S_2$ in the obvious way.\\

Step 1. $\Sc^*(Z^n(F))^{\widetilde{GL_n \times GL_n}(F),\chi}=0$.\\
Consider the projection on the first coordinate from $Z^n$ to the
transitive $\widetilde{GL_n \times GL_n}$-space  $Q^n=GL_n$.
Choose the point $Id \in Q^n$. Its stabilizer is $\tH$ and its
fiber is $W$. Hence by Frobenius descent (Theorem \ref{Frob}),
$\Sc^*(Z^n(F))^{\widetilde{GL_n(F) \times GL_n(F)},\chi} \cong
\Sc^*(W(F))^{\tH(F),\chi}$ which is zero by the key lemma.\\

Step 2. $\Sc^*(Z^n(F))^{\tG(F),\chi}=0$.\\
Consider the space $Y:=Z^n \times F^{\times}$ and let the group
$GL_n \times GL_n$ act on it by $(h_1,h_2) (z,\lambda):=
((h_1,h_2)z,\det h_1 \det h_2^{-1} \lambda)$. Extend this action
to action of $\widetilde{GL_n \times GL_n}$ by $\sigma(z,\lambda)
:= (\sigma(z),\lambda)$. Consider the projection $Z^n \times
F^{\times} \to F^{\times}$. By Frobenius descent (Theorem
\ref{Frob}),
$$\Sc^*(Y(F))^{\widetilde{GL_n(F) \times GL_n(F)},\chi} \cong
\Sc^*(Z^n(F))^{\tG,\chi}.$$ Let $Y'$ be equal to $Y$ as an
$l$-space and let $\widetilde{GL_n \times GL_n}$ act on $Y'$ by
$(h_1,h_2) (z,\lambda):= ((h_1,h_2)z,\lambda)$ and
$\sigma(z,\lambda):=(\sigma(z),\lambda)$. Now $Y$ is isomorphic to
$Y'$ as a $\widetilde{GL_n \times GL_n}$ space by
$((A,v,\phi),\lambda) \mapsto ((A,v,\phi), \lambda \det A^{-1})$.

Since $\Sc^*(Z^n(F))^{\widetilde{GL_n(F) \times GL_n(F)},\chi} =
0$, Proposition \ref{Product} implies that
$\Sc^*(Y'(F))^{\widetilde{GL_n(F) \times GL_n(F)},\chi} = 0$ and
hence $\Sc^*(Y(F))^{\widetilde{GL_n(F) \times GL_n(F)},\chi} = 0$
and thus
$\Sc^*(Z^n(F))^{\tG_n(F),\chi}=0$.%
\end{proof}
Till the end of the subsection we assume that $F$ is Archimedean,
since the argument in the non-Archimedean case is similar but much
easier.

\begin{corollary} \label{SymW}
$\Sc^*(W_n(F), Sym^k(\mathrm{gl}_n^*)(F))^{\tG(F),\chi}=0.$
\end{corollary}
\begin{proof}
Consider the Killing form $K:\mathrm{gl}_n^* \to \mathrm{gl}_n$.
Let $U:=K^{-1}(Q^n_n)$. In the same way as in the previous
corollary one can show that $\Sc^*(W_n(F) \times
U(F))^{\tG(F),\chi}=0$. Hence by Proposition \ref{Trick},
$\Sc^*(W_n(F), Sym^k(\mathrm{gl}_n^*)(F))^{\tG(F),\chi}=0.$
\end{proof}

\begin{corollary}
We have
\begin{multline*}\Sc^*(W_i(F) \times W_{n-i}(F),Sym^k(0
\times \mathrm{gl}_{n-i}^*)(F))^{H_i(F) \times G_{n-i}(F)} =\\=
\Sc^*(W_i(F) \times W_{n-i}(F),Sym^k(0 \times
\mathrm{gl}_{n-i}^*)(F))^{\tH_i(F) \times \tG_{n-i}(F)}.
\end{multline*}

\end{corollary}
\begin{proof}
It follows from Lemma \ref{keyGLWeak}, the last corollary and
Proposition \ref{Product}.
\end{proof}

Now we are ready to prove Proposition \ref{Enough}.
\begin{proof}[Proof of Proposition \ref{Enough}]
Fix $i<n$. Consider the projection $pr_1:Z^i \to Q^i$. It is easy
to see that the action of $\tG$ on $Q^i$ is transitive. Denote
$$A_i := \left(
  \begin{array}{cc}
    Id_{i\times i} & 0 \\
    0 & 0 \\
  \end{array}
\right) \in Q^i.$$ Denote by $G_{A_i}:= \mathrm{Stab}_G(A_i)$ and
$\tG_{A_i}:= \mathrm{Stab}_{\tG}(A_i)$. Note that they are
unimodular. By Frobenius descent (Theorem \ref{Frob}),
$$\Sc^*(Z^i(F),Sym^k(CN_{Z^i}^X
)(F))^{\tG(F),\chi}=\Sc^*(W(F),Sym^k(CN_{Q^i,A_i}^{\mathrm{gl}_n})(F))^{\tG_{A_i}(F),\chi}.$$
Hence it is enough to show that
$$\Sc^*(W(F),Sym^k(CN_{Q^i,A_i}^{\mathrm{gl}_n})(F))^{G_{A_i}(F)}=\Sc^*(W(F),Sym^k(CN_{Q^i,A_i}^{\mathrm{gl}_n})(F))^{\tG_{A_i}(F)}.$$

It is easy to check by explicit computation that
\begin{itemize}
\item
$H_i \times G_{n-i}$ is canonically
embedded into $G_{A_i}$, 
\item
$W$ is isomorphic to $W_i \times W_{n-i}$ as $H_i \times
G_{n-i}$-spaces
\item
$CN_{Q^i,A_i}^{\mathrm{gl}_n}$ is isomorphic to $0 \times
\mathrm{gl}_{n-i}^*$ as $H_i \times G_{n-i}$ representations.
\end{itemize}
Let $\xi \in
\Sc^*(W,Sym^k(CN_{Q^i,A_i}^{\mathrm{gl}_n})(F))^{G_{A_i}(F)}$. By
the previous corollary, $\xi$ is $\tH_i(F) \times
\tG_{n-i}(F)$-invariant. Since $\xi$ is also
$G_{A_i}(F)$-invariant, it is $\tG_{A_i}(F)$-invariant.
\end{proof}
\subsection{Proof of Lemma \ref{keyGLWeak}} \label{ProofKeyGLWeak}

\begin{proposition} \label{1Enough}
It is enough to prove the key lemma for $n=1$.
\end{proposition}
\begin{proof}
Consider the subgroup $T_n \subset H_n$ consisting of diagonal
matrices, and $\widetilde{T}_n:=T_n\rtimes S_2 \subset \tH_n$. It
is enough  to prove $\Sc^*(W_n(F))^{\widetilde{T}_n(F),\chi}=0.$

Now by Proposition \ref{Product} it is enough to prove
$\Sc^*(W_1(F))^{\tH_1(F),\chi}=0$.
\end{proof}

From now on we fix $n:=1$, $H:=H_1$, $\tH:=\tH_1$ and $W:=W_1$.
Note that $H(F) = F^{\times}$ and $W(F)=F^2$. The action of $H(F)$
is given by $\rho(\lambda)(x,y):=(\lambda x, \lambda^{-1} y)$ and
extended to the action of $\tH(F)$ by the involution $\sigma(x,y)
= (y,x)$.

Let $Y:=\{(x,y)\in  F^2| \, xy=0\} \subset W$ be the {\bf cross}
and $Y':=Y \setminus \{0\}$.
\begin{lemma}\label{AllOnCross}
Every $(\tH(F),\chi)$-equivariant distribution on $W$ is supported
inside the cross $Y$.
\end{lemma}

\begin{proof}
Denote $U:=W \setminus Y$. We have to show
$\Sc^*(U(F))^{\tH(F),\chi}=0$. Consider the coordinate change
$U(F) \cong F^{\times} \times F^{\times}$ given by $(x,y) \mapsto
(xy,x/y)$. It is an isomorphism of $\tH(F)$-spaces where the
action of $\tH(F)$ on $F^{\times} \times F^{\times}$ is only on
the second coordinate, and given by $\lambda (w) = \lambda^2 w$
and $\sigma(w) = w^{-1}$. Clearly,
$\Sc^*(F^\times)^{\tH(F),\chi}=0$ and hence by Proposition
\ref{Product} $\Sc^*(F^\times \times F^\times)^{\tH(F),\chi}=0$.
\end{proof}

%

\begin{lemma} $ $\\
(i) $\Sc^*(W(F)\setminus Y(F))^{\tH(F),\chi}=0$.\\
(ii) Any distribution $\xi \in \Sc^*(Y'(F))^{\tH(F),\chi}$ is
invariant with respect to homotheties.

Moreover, if $F$ is Archimedean then for all $k \in \Z_{\geq 0}$,
any distribution $\xi \in \Sc^*(Y'(F),
Sym^k(CN_{Y'}^W)(F))^{\tH(F),\chi}$ is $\R$-homogeneous
of type $\alpha_k$ where $\alpha_k(\lambda):=\lambda^{-2k}$.\\
(iii) $\Sc^*(\{0\})^{\tH(F),\chi}=0$. Moreover, if $F$ is
Archimedean then
$\Sc^*(\{0\},Sym^k(CN_{\{0\}}^W)(F))^{\tH(F),\chi}=0$.
\end{lemma}
If $F$ is non-Archimedean this proposition is clear. Let us prove
it for Archimedean $F$.
\begin{proof}
We have proven (i) in the proof of the previous lemma.\\
(ii) Fix $x_0:=(1,0)\in Y'(F)$. Now we want to use Frobenius
descent (Theorem \ref{Frob}) for the group $\tH(F) \times
\R^{\times}$ and character $\chi \times \alpha_k$. Note that
$Stab_{\tH}(x_0)$ is trivial and $Stab_{\tH(F) \times
\R^{\times}}(x_0) \cong \R^{\times}$. Note that $N_{Y',x_0}^W(F)
\cong F$ and $Stab_{\tH(F) \times \R^{\times}}(x_0)$ acts on it by
$\rho(\lambda)a = \lambda^{2} a$. So we have $$
Sym^k(N_{Y',x_0}^W(F))=Sym^k(N_{Y',x_0}^W(F))^{\R^{\times},\alpha_k^{-1}}.$$
So by Frobenius descent any distribution $\xi \in \Sc^*(Y'(F),
Sym^k(CN_{Y'}^W)(F))^{\tH(F),\chi}$ is $\R$-homogeneous of type
$\alpha_k$.\\
(iii) is a simple computation. Also, it can be deduced from (i)
using Proposition \ref{Trick}.
\end{proof}

\begin{proof}[Proof of Lemma \ref{keyGLWeak}]
Let $\xi \in \Sc^*(W(F))^{\tH(F),\chi}$. Consider the quadratic
form on $W=F^2$ given by $B(x,y):=xy$. It defines a Fourier
transform on $\Sc^*(W(F))$. By Lemma \ref{FouInv}, $\Fou(\xi) \in
\Sc^*(W(F))^{\tH(F),\chi}$.

Hence by Lemma \ref{AllOnCross}, $\Supp \xi, \Supp \Fou(\xi)
\subset Y$. By Homogeneity Theorem (Theorem \ref{Homog}) this
implies that $\xi$ is $B$-adapted.

Finally, the previous lemma and Theorem \ref{NashFilt} imply that
all $B$-adapted distributions in $\Sc^*(W(F))^{\tH(F),\chi}$ are
0.
\end{proof}

\section{The pair $(\mathrm{O}_{n+1}(F),\mathrm{O}_n(F))$}
Let $(W,Q)$ be a quadratic space defined over $F$ and fix $e \in
W$ a unit vector. Consider the quadratic space $V={e}^{\bot}$ with
$q=Q|_{V}$. Define the standard imbedding $\mathrm{O}(V)
\hookrightarrow \mathrm{O}(W)$ and consider the two-sided action
of $\mathrm{O}(V) \times \mathrm{O}(V)$ on $\mathrm{O}(W)$ defined
by
$(g_1,g_2)h:=g_1 hg_2^{-1}$. We also consider the anti-involution
$\sigma$ of $O_{Q}$ given by $\sigma(g)=g^{-1}$. In this paper we
prove the following Theorem \begin{theorem} \label{OWeak} Any
$\mathrm{O}(V) \times \mathrm{O}(V)$ invariant distribution on
$\mathrm{O}(W)$ is invariant under $\sigma$.
\end{theorem}

By Theorem \ref{DistCrit} and Corollary \ref{GKCorWeakPairs}  this
theorem implies

\begin{theorem}
The pair $(\mathrm{O}(W),\mathrm{O}(V))$ is a Gelfand pair.
\end{theorem}

\begin{remark}
In fact, our proof can be easily adapted to prove also an
analogous theorem for unitary groups.
\end{remark}

\begin{remark}
This result was earlier proven in a bit different form in
\cite{vD} for $F=\R$, in \cite{Apa-vD} for $F=\C$ and in
\cite{Bos-vD} for non-Archimedean $F$ of characteristic zero. Here
we give a different proof.

Also, an analogous theorem for unitary groups over $\R$ is proven
in \cite{vD2}.
\end{remark}

For the proof we need some further notations.

\begin{itemize}
\item
$O_{q}=O(V,q)$ is the group of isometries of the quadratic space
$(V,q)$.

\item
 $G_{q}=O(V,q) \times O(V,q)$.
\item
 $\Delta: O_{q} \to
G_{q}$ the diagonal. $H_{q}=\Delta(O_{q}) \subset G_{q}$.

\item
$\sigma(g_{1},g_{2})=(g_{2},g_{1})$.

\item
$\widetilde{G_{q}}=G_{q} \rtimes \{1,\sigma\}$, same for
$\widetilde{H_{q}}$

\item
$\chi : \widetilde{G_{q}} \to \{+1,-1\}$ the non trivial character
with $\chi(G_{q})=1$.

\item
$\widetilde{G_{Q}}$ acts on $O_{Q}$ by
$(g_{1},g_{2})x=g_{1}xg_{2}^{-1}$ and $\sigma(x)=x^{-1}$.

\end{itemize}

By Theorem \ref{NoSNoDist}, Theorem \ref{OWeak} follows from the
following theorem:

\begin{theorem} \label{OReform}
$$\Sc^*(O_{Q})^{\widetilde{G_{q}},\chi}=0.$$
\end{theorem}

\subsection{Proof of Theorem \ref{OReform}}
We denote by $\Gamma=\{w \in W: Q(w)=1\}$. Note that by Witt's
theorem $\Gamma$ is an $O_{Q}$ transitive set and therefore
$\Gamma \times \Gamma$ is a transitive $\widetilde{G}_{Q}$ set
where the action of $G_{Q}$ is the standard action on $W \oplus W$
and $\sigma$ acts by flip.

Applying Frobenuis descent (Theorem \ref{Frob}) to projections of
$O_Q \times \Gamma \times \Gamma$ first on $\Gamma \times \Gamma$
and then on $O_Q$ we have
$$\Sc^*
(O_Q)^{\widetilde{G_{q}},\chi}=\Sc^*(O_{Q} \times \Gamma \times
\Gamma)^{\widetilde{G_{Q}},\chi}$$ and also that

$$\Sc^*(O_{Q} \times \Gamma
\times \Gamma)^{\widetilde{G_{Q}},\chi}=\Sc^*(\Gamma \times
\Gamma)^{\widetilde{H_{Q}},\chi}$$

In what follows we will abuse notation and write $Q(u,v)$ for the
bilinear form defined by $Q$. Define a map $D: \Gamma \times
\Gamma \to Z$ where $Z=\{(v,u) \in W \oplus W: Q(v,u)=0 ,
Q(v+u)=4\}$ by
$$D(x,y)=(x+y,x-y).$$ $D$ defines an $\widetilde{G}_{Q}$-equivariant
homeomorphism and thus we need to show that
$$\Sc^*(Z)^{\widetilde{H_{Q}},\chi}=0$$
Here, the action of $G_{Q}$ on $Z \subset W \oplus W$ is the
restriction of its action on $W \oplus W$ while the action of
$\sigma$ is given by $\sigma(v,u)=(v,-u).$

Now we cover $Z = U_{1} \cup U_{2}$ where
$$U_{1}=\{(v,u) \in Z: Q(v) \ne 0\}$$ and $$U_{2}=\{(v,u) \in Z: Q(u) \ne 0\}$$
We will show $\Sc^*(U_{1})^{\widetilde{H_{Q}},\chi}=0$, and the
proof for $U_{2}$ is analogous. This will finish the proof.

\begin{lemma}\label{non uniform}
$\Sc^*(U_{1})^{\widetilde{H_{Q}},\chi}=0$ \end{lemma}

\begin{proof}[Proof for non-archimedean $F$.]

Consider $\ell_{1}:U_{1} \to F-\{0\}$ defined as
$\ell_{1}(v,u)=Q(v)$. By the Localization Principle (Corollary
\ref{LocPrinSub}), it is enough to show
$\Sc^*(U_{1}^{\alpha})^{\widetilde{H_{Q}},\chi}=0$ where
$U_{1}^{\alpha}=\ell_{1}^{-1}(\alpha)$, for any $\alpha \in
F-\{0\}$. But

$$U_{1}^{\alpha}=\{(v,u)| Q(v)=\alpha, Q(u)=4-\alpha, Q(v,u)=0\}$$

Let $W^{\alpha}=\{w \in W| Q(w)=\alpha\}$ and let
$p_1:U_{1}^{\alpha} \to W^{\alpha}$ be given by $p_1(v,u)=v$.

On $W^{\alpha}$ our group acts transitively. Fix a vector $v_{0}
\in W^{\alpha}$.\\ Denote $H(v_0):= H_{(Q|_{v_{0}^{\bot}})}$ and
$\widetilde{H}(v_0):= \widetilde{H}_{(Q|_{v_{0}^{\bot}})}$.

The stabilizer in $\widetilde{H}_{Q}$ of $v_{0}$ is
$\widetilde{H}(v_0)$. The fiber $p_1^{-1}(v_{0})=\{a \in
{v_{0}^{\bot}}| Q(a)=4-\alpha\}.$

Frobenius descent implies that

$$\Sc^*(U_{1}^{\alpha})^{\widetilde{H_{Q}},\chi}=\Sc^*(p_1^{-1}(v_{0}))^{\widetilde{H}(v_0),\chi}$$
But clearly $\Sc^*(p_1^{-1}(v_{0}))^{\widetilde{H}(v_0),\chi}=0$
as $-Id \in H(v_0)$.
\end{proof}

\begin{proof}[Proof for archimedean $F$.]
Now let us consider the archimedean case. Define $U:=\{(v,u) \in
U_1 | u \neq 0 \}$. Note that the map $\ell _1|_U$ is a
submersion, so the same argument as in the non-archimedean case
shows that $\Sc^*(U)^{\widetilde{H_{Q}},\chi}=0$. Let $Y:=\{(v\in
W | Q(v) = 4\} \times \{0\}$ be the complement to $U$ in $U_1$.
By Theorem \ref{NashFilt}, it is enough to prove
$$\Sc^*(Y,Sym^k(CN_Y^{U_1}))^{\widetilde{H_{Q}},\chi}=0.$$

Note that the action of $\widetilde{H}_{Q}$ on $Y$ is transitive,
and fix a point $(v,0) \in Y$. The stabilizer in
$\widetilde{H}_{Q}$ of $(v,0)$ is $\widetilde{H}(v)$, and the
normal space to $Y$ at $(v,0)$ is $v^{\bot}$. So Frobenius descent
(Theorem \ref{Frob}) implies that
$$\Sc^*(Y,Sym^k(CN_Y^{U_1}))^{\widetilde{H_{Q}},\chi}=Sym^k(v^{\bot})^{\widetilde{H}(v),\chi}.$$
But clearly $Sym^k(v^{\bot})^{\widetilde{H}(v),\chi}=0$ as $-Id
\in H(v)$.
\end{proof}

\chapter{Proof of the main results}\label{ChProofMain}

In this chapter we assume that $F$ has zero characteristic. This
chapter is based on \cite{AG_AMOT,Aiz}.

Consider the standard imbedding $\mathrm{GL}_n(F) \hookrightarrow
\mathrm{GL}_{n+1}(F)$. We consider the action of
$\mathrm{GL}_n(F)$ on $\mathrm{GL}_{n+1}(F)$ by conjugation. In
this chapter we prove Theorem \ref{StrongDist}, namely
\begin{theorem*}
Any $\mathrm{GL}_n(F)$ - invariant distribution on
$\mathrm{GL}_{n+1}(F)$ is invariant with respect to transposition.
\end{theorem*}

By Corollary \ref{StrongDist2StrongRep} it implies Theorem
\ref{StrongRep}, namely
\begin{theorem*}
$(\mathrm{GL}_{n+1}(F),\mathrm{GL}_n(F))$ is a strong Gelfand
pair.
\end{theorem*}

\section{Structure of the proof}
We will now briefly sketch the main ingredients of our proof of
Theorem \ref{StrongDist}.

First we show that we can switch to the following problem. The
group $\mathrm{GL}_n(F)$ acts on a certain linear space $X_n$ and
$\sigma$ is an involution of $X_n$. We have to prove that every
$\mathrm{GL}_n(F)$-invariant distribution on $X_n$ is also
$\sigma$-invariant. We do that by induction on $n$. Using the
Harish-Chandra descent method we show that the induction
hypothesis implies that this holds for distributions on the
complement to a certain small closed subset $S \subset X_n$. We
call this set {\it the singular set}.

Next we assume the contrary: there exists a non-zero
$\mathrm{GL}_n(F)$-invariant distribution $\xi$ on $X$ which is
anti-invariant with respect to $\sigma$.

We use the notion of singular support of a distribution (see
section \ref{SingSupp}) . Let $T \subset T^*X$ denote the singular
support of $\xi$. Using Fourier transform and the fact any such
distribution is supported in $S$ we obtain that $T$ is contained
in $\check{S}$ where $\check{S}$ is a certain small subset in
$T^*X$.

Then we use a powerful tool, Theorem \ref{ThmSingSupp}, which
states that the singular support of a distribution is a weakly
coisotropic variety in the cotangent bundle. This enables us to
show, using a complicated but purely geometric argument, that the
support of $\xi$ is contained in a much smaller subset of $S$.

Finally it remains to prove that any $\mathrm{GL}_n(F)$-invariant
distribution that is supported on this subset together with its
Fourier transform is zero. This is proven using Homogeneity
Theorem (Theorem \ref{ArchHom}) which in turn uses Weil
representation.\\

Now let us describe the chapter section by section.

%

In section \ref{SpecNat} we introduce notation that we will use
in our proof.\\

In section \ref{HC} we use the Harish-Chandra descent method.

In subsection \ref{Lin} we linearize the problem to a problem on
the linear space $X=\sl(V)\times V \times V^*$, where $V=F^n$.

In subsection \ref{subHC} we perform the Harish-Chandra descent on
the $\sl(V)$-coordinate and $V\times V^*$ coordinate separately
and then use non-linear automorphisms $\nu_{\lambda}$ of $X$ to
descend further to the singular set $S$. The automorphisms
$\nu_{\lambda}$ of $X$ first appeared in the first proof of
Theorem \ref{StrongDist} for non-Archimedean fields, in
\cite{AGRS}, and were essential part of the proof. In the proof we
give here they are used mainly to shorten computations. However,
in the proof of an analogous theorem for orthogonal groups over
Archimedean fields in \cite{SZ} they are again used in
an essential way.\\

In section \ref{Red2Geo} we reduce Theorem \ref{StrongDist} to the
following geometric statement: any coisotropic subvariety of
$\check{S}$ is contained in a certain set $\check{C}_{X \times
X}$. The reduction is done using the fact that the singular
support of a distribution has to be coisotropic, and the following
proposition: any $\mathrm{GL}(V)$-invariant distribution on $X$
such that it and its Fourier transform are supported on
$\sl(V)\times (V \times 0 \cup 0 \times V^*)$ is zero.

In subsection \ref{ProofLemCros} we prove this proposition using
Homogeneity theorem.\\

In section \ref{ProofGeo} we prove the geometric statement. On one
hand, this section is rather complicated. On the other hand, it is
purely geometric and involves no distributions. Using the fact
that the only non-distinguished orbit in $sl_n$ is the regular
orbit, we reduce the geometric statement to the Key Lemma that
states that a certain subset $R_A$ of $V \times V^* \times V
\times V^*$ contains no non-empty coisotropic subvarieties. This
subset has Lagrangian dimension, i.e. $\dim R_A =2n$. Hence in
order to prove the Key Lemma it is enough to exhibit one
additional equation on $R_A$. We write this equation explicitly.

%
%


\section{Notation} \label{SpecNat}  $ $\\
In this chapter we will use the following notation
\itemize{
\item Let $V:=V_n$ be the standard $n$-dimensional linear space defined over $F$.

\item Let $\sl(V)$ denote the Lie algebra of operators with zero
trace. \item Denote $X:=X_n:=\sl(V_n) \times V_n \times V_n^*$
\item $G:=G_n:=\mathrm{GL}(V_n)$ \item $\g:=\g_n:=\Lie
(G_n)=\mathrm{gl}(V_n)$ \item $\tG:= \tG_n:= G_n \rtimes
\{1,\sigma\}$, where the action of the 2-element group
$\{1,\sigma\}$ on $G$ is given by the involution $g \mapsto
{g^t}^{-1}$. \item We define a character $\chi$ of $\tG$ by
$\chi(G)= \{1 \}$ and $\chi(\tG - G)= \{-1\}$.
\item Let $G_n$ act on $G_{n+1}$, $\g_{n+1}$ and on $\sl(V_n)$ by
$g(A):= gAg^{-1}$.
\item Let  $G$ act on $V \times V^*$ by $g(v,\phi):=(gv,
({g}^{*})^{-1} \phi)$. This gives rise to an action of $G$ on $X$.
\item Extend the actions of $G$ to actions of $\tG$ by
$\sigma(A):=A^t$ and $\sigma(v,\phi):=(\phi^t,v^t)$. \item We
consider the standard scalar products on $\sl(V)$ and $V \times
V^*$. They give rise to a scalar product on $X$. \item We identify
the cotangent bundle $T^*X$ with $X \times X$ using the above
scalar product.
\item Let $\cN:=\cN_n \subset \sl(V_n)$ denote the cone of nilpotent operators.
\item $C := (V \times 0) \cup (0 \times V^* ) \subset V \times
V^*$.

\item $\check{C} := (V \times 0  \times V \times 0) \cup (0 \times
V^* \times 0 \times V^*) \subset V \times V^* \times V \times
V^*$.

\item $\check{C}_{X \times X} := (\sl(V) \times V \times 0 \times
\sl(V) \times V \times 0) \cup (\sl(V) \times 0 \times V^* \times
\sl(V) \times 0 \times V^*) \subset X \times X$.

\item $S:= \{(A,v,\phi) \in X_n | A^n=0 \text{ and } \phi(A^i v)=0 \text{ for any } 0 \leq i \leq n\}$.
\item \begin{multline*}
\check{S} := \{((A_1,v_1,\phi_1),(A_2,v_2,\phi_2)) \in X
\times X \, | \, \forall i,j \in \{1,2\} \\
(A_i,v_j,\phi_j) \in S \text{ and } \forall \alpha \in
\mathrm{gl}(V), \alpha (A_1,v_1,\phi_1) \bot (A_2,v_2,\phi_2)\}
\end{multline*}.
\item Note that
\begin{multline*}
\check{S} = \{((A_1,v_1,\phi_1),(A_2,v_2,\phi_2)) \in X \times X \, | \, \forall i,j \in \{1,2\} \\
(A_i,v_j,\phi_j) \in S \text{ and } [A_1,A_2] + v_1 \otimes \phi_2
- v_2 \otimes \phi_1 =0\}.
\end{multline*}
\item $\check{S}':= \check{S} - \check{C}_{X \times X}$.
\item $\Gamma := \{ (v,\phi) \in V \times V^* \, | \, \phi(V)=0 \}$.
\item
For any $\lambda \in F$ we define $\nu_\lambda:X \to X$ by
$\nu_{\lambda}(A,v,\phi):=(A+\lambda v\otimes \phi-\lambda
\frac{\langle \phi,v \rangle}{n}\mathrm{Id},v,\phi).$
\item It defines $\check{\nu}_\lambda : X \times X \to X \times X$. It is given by
\begin{multline*}
\check{\nu}_\lambda ((A_1,v_1,\phi_1),(A_2,v_2,\phi_2)) = \\=
((A_1+\lambda v_1\otimes \phi_1-\lambda \frac{\langle \phi_1,v_1
\rangle}{n}\mathrm{Id},v_1,\phi_1), (A_2, v_2 -\lambda A_2
v_1,\phi_2 - \lambda A_2^* \phi_1)).
\end{multline*}
}

%
%

\section{Harish-Chandra descent} \label{HC}


%

\subsection{Linearization} \label{Lin}$ $\\
In this subsection we reduce Theorem A to the following one
\begin{theorem}  \label{descendant}
$\Sc^*(X(F))^{\tG(F),\chi}=0.$
\end{theorem}

We will divide this reduction to several propositions.\\

\begin{proposition} \label{Red1}
If  $\cD(G_{n+1}(F))^{\tG_n(F),\chi}=0$ then Theorem A holds.
\end{proposition}
The proof is straightforward.

\begin{proposition} 
If  $\Sc^*(G_{n+1}(F))^{\tG_n(F),\chi}=0$ then
$\cD(G_{n+1}(F))^{\tG_n(F),\chi}=0$.
\end{proposition}
Follows from Theorem \ref{NoSNoDist}.

\begin{proposition} \label{Red2}
If  $\Sc^*(\g_{n+1}(F))^{\tG_n(F),\chi}=0$ then
$\Sc^*(G_{n+1}(F))^{\tG_n(F),\chi}=0$.
\end{proposition}
\begin{proof} Let $\xi \in \Sc^*(G_{n+1}(F))^{\tG_n(F),\chi}$. We have to
prove $\xi=0$. Assume the contrary. Take $p \in
\mathrm{Supp}(\xi)$. Let $t=\mathrm{det}(p)$. Let $f\in \Sc(F)$ be
such that $f$ vanishes in a neighborhood of 0 and $f(t) \neq 0$.
Consider the determinant map $\mathrm{det}:G_{n+1}(F) \to F$.
Consider $\xi':=(f \circ \mathrm{det})\cdot \xi$. It is easy to
check that $\xi' \in \Sc^*(G_{n+1}(F))^{\tG_n(F),\chi}$ and $p \in
\mathrm{Supp}(\xi')$. However, we can extend $\xi'$ by zero to
$\xi'' \in \Sc^*(\g_{n+1}(F))^{\tG_n(F),\chi}$, which is zero by
the assumption. Hence $\xi'$ is also zero. Contradiction.
\end{proof}

\begin{proposition}
If $\Sc^*(X_{n}(F))^{\tG_n(F),\chi}=0$ then
$\Sc^*(\g_{n+1}(F))^{\tG_n(F),\chi}=0$.
\end{proposition}

\begin{proof}
The $\tG_n(F)$-space $\mathrm{gl}_{n+1}(F)$ is isomorphic to
$X_n(F) \times F \times F$ with trivial action on $F\times F$.
This isomorphism is given by
$$ \left(
  \begin{array}{cc}
    A_{n \times n} & v_{n\times 1} \\
    \phi_{1\times n} & \lambda \\
  \end{array}
\right) \mapsto ((A - \frac{\tr A}{n} \mathrm{Id}, v
,\phi),\lambda, \tr A).$$
 \end{proof}

\subsection{Harish-Chandra descent} \label{subHC}$ $\\

Now we start to prove Theorem \ref{descendant}. The proof is by
induction on $n$. Till the end of the paper we will assume that
Theorem \ref{descendant} holds for all $k<n$ for both archimedean
local fields.

The theorem obviously holds for $n=0$. Thus from now on we assume
$n\geq 1$. The goal of this subsection is to prove the following
theorem.

\begin{theorem}  \label{LastRed}
$\Sc^*(X(F) - S(F))^{\tG(F),\chi}=0.$
\end{theorem}

In fact, one can prove this theorem directly using  Theorem
\ref{HC_Thm}. However, this will require long computations. Thus,
we will
 divide the proof to several steps and use some
tricks to avoid part of those computations.

\begin{proposition}  
$\Sc^*(X(F) - (\cN \times V \times V^*)(F))^{\tG(F),\chi}=0.$
\end{proposition}
\begin{proof}

By Theorem \ref{HC_Thm} it is enough to prove that for any
semisimple $A \in \sl(V)$ we have
$$\Sc^*((N_{GA,A}^{\sl(V)} \times (V \times V^*)) (F))^{\tG(F)_A,\chi} = 0.$$

Now note that $\tG(F)_A \cong \prod \tG_{n_i}(F_i)$ where $n_i<n$
and $F_i$ are some field extensions of $F$.
Note also that
$$(N_{GA,A}^{\sl(V)} \times V \times V^*)(F) \cong
\sl(V)_{A} \times (V \times V^*)(F) \cong \prod X_{n_i}(F_i)
\times {\mathcal Z}(\sl(V)_{A})(F),$$ where ${\mathcal
Z}(\sl(V)_{A})$ is the center of $\sl(V)_{A}$. Clearly, $\tG_A$
acts trivially on ${\mathcal Z}(\sl(V)_{A})$.

Now by Proposition \ref{Product} the induction hypothesis implies
that
$$\Sc^*(\prod
X_{n_i}(F_i) \times {\mathcal Z}(\sl(V)_{A})(F))^{\prod
\tG_{n_i}(F_i),\chi} = 0.$$
\end{proof}

In the same way we obtain the following proposition.

\begin{proposition}  
$\Sc^*(X(F) - (\sl(V) \times \Gamma)(F))^{\tG(F),\chi}=0.$
\end{proposition}

\begin{corollary}  \label{PartDes}
$\Sc^*(X(F) - (\cN \times \Gamma)(F))^{\tG(F),\chi}=0.$
\end{corollary}

%
\begin{lemma}
Let $A \in \sl(V)$, $v \in V$ and $\phi \in V^*$. Suppose $A +
\lambda v \otimes \phi$ is nilpotent for all $\lambda \in F$. Then
$\phi(A^i v)=0$ for any $i \geq 0$.
\end{lemma}
\begin{proof}
Since $A + \lambda v \otimes \phi$ is nilpotent, we have $tr(A+
\lambda v \otimes \phi)^k=0$ for any $k \geq 0$ and $\lambda \in
F$. By induction on $i$ this implies that $\phi(A^i v)=0$.
\end{proof}

\begin{proof}[Proof of Theorem \ref{LastRed}]
By the previous lemma, $\bigcap_{\lambda \in F}\nu_{\lambda}(\cN
\times \Gamma) \subset S$. Hence $\bigcup _{\lambda \in
F}\nu_{\lambda}(X - \cN \times \Gamma) \supset X-S$.

By Corollary \ref{PartDes} $\Sc^*(X(F) - (\cN \times
\Gamma)(F))^{\tG(F),\chi}=0$. Note that $\nu_{\lambda}$ commutes
with the action of $\tG$. Thus $\Sc^*(\nu_{\lambda}(X(F) - (\cN
\times \Gamma)(F)))^{\tG(F),\chi}=0$ and hence $\Sc^*(X(F) -
S(F))^{\tG(F),\chi}=0.$
\end{proof}

\begin{remark}
The automorphisms $\nu_{\lambda}$ of $X$ first appeared in the
first proof of Theorem \ref{StrongDist} for non-Archimedean
fields, in \cite{AGRS}, and were essential part of the proof. In
the proof we give here they are used mainly to shorten
computations. However, in the proof of an analogous theorem for
orthogonal groups over Archimedean fields in \cite{SZ} they are
again used in an essential way.
\end{remark}

\begin{remark}
The main ingredient of $\nu_{\lambda}$, i.e. the map $(A,v,\phi)
\mapsto A+v\otimes \phi$ is the moment map corresponding to the
action of $\G$ on $\g \times V \times V^*$.
\end{remark}

\section{Reduction to the geometric statement} \label{Red2Geo}

In this section coisotropic variety means $X \times X$-coisotropic
variety.

The goal of this section  is to reduce Theorem \ref{descendant} to
the following geometric statement.
\begin{theorem}[geometric statement] \label{GeoStat}
For any coisotropic subvariety of $T \subset \check{S}$ we have
$T\subset \check{C}_{X \times X}$.
\end{theorem}

Till the end of this section we will assume the geometric
statement.

\begin{proposition}
Let $\xi \in \Sc^*(X(F))^{\tG(F),\chi}=0.$ Then $\Supp(\xi)
\subset (\sl(V) \times C)(F)$.
\end{proposition}
\begin{proof}[Proof for the case $F=\R$] $ $

Step 1. $SS(\xi) \subset \check{S}$.\\
We know that $$\Supp(\xi), \Supp(\Fou^{-1}_{\sl(V)}\xi),
\Supp(\Fou^{-1}_{V \times V^*}(\xi)), \Supp(\Fou^{-1}_{X}(\xi))
\subset S(F).$$
By property (\ref{Fou}) of the singular support this implies that
$$SS(\xi) \subset (S \times X) \cap F_{\sl(V)}(S \times X) \cap
F_{V\times V^*}(S \times X) \cap F_{X}(S \times X).$$ On the other
hand, $\xi$ is $G(F)$-invariant and hence by property (\ref{Ginv})
of the singular support
$$SS(\xi) \subset \{((x_1,x_2) \in X \times X \, | \, \forall g
\in \g, g(x_1) \bot x_2 \}. $$ Thus $SS(\xi)\subset \check{S}$.

Step 2. $SS(\xi) \subset \check{C}_{X \times X}$.\\
By Corollary \ref{Gaber}, $SS(\xi)$ is $X \times X$-coisotropic
and hence by the geometric statement $SS(\xi) \subset \check{C}_{X
\times X}$.

Step 3. $\Supp(\xi) \subset (\sl(V) \times C)(F)$.\\
Follows from the previous step by property (\ref{Supp2SS}) of the
singular support.
\end{proof}

The case $F=\C$ is proven in the same way using the following
corollary of the geometric statement.
\begin{proposition}
Any $(X \times X)_{\C}$-coisotropic subvariety of $\check{S}_{\C}$
is contained in $(\check{C}_{X \times X})_{\C}$.
\end{proposition}

Now it is left to prove the following proposition.

\begin{proposition} \label{Cross}
Let $\xi \in \Sc^*(X(F))^{\tG(F),\chi}$ be such that
$$\Supp(\xi),\Supp(\Fou_{V \times V^*}(\xi)) \subset (\sl(V) \times C)(F).$$
Then $\xi=0$.
\end{proposition}

\subsection{Proof of Proposition \ref{Cross}} \label{ProofLemCros} $ $\\
Proposition \ref{Cross} follows from the following lemma.

\begin{lemma}
Let $F^{\times}$ act on $V\times V^*$ by $\lambda(v,\phi) :=
(\lambda v, \frac{\phi}{\lambda})$.
 Let $\xi \in
\Sc^*((V\times V^*)(F))^{F^{\times}}$ be such that $$\Supp(\xi),
\Supp(\Fou_{V \times V^*}(\xi)) \subset C(F).$$ Then $\xi=0$.
\end{lemma}

By Homogeneity Theorem (Theorem \ref{ArchHom}) it is enough to
prove the following lemma.

\begin{lemma}
Let $\mu$ be a character of $F^{\times}$ given by $|| \cdot ||^n
u$ or $|| \cdot ||^{n+1} u$ where $u$ is some unitary character.
Let $F^{\times} \times F^{\times}$ act on $V \times V^*$ by
$(x,y)(v,\phi) = (\frac{y}{x} v, \frac{1}{xy} \phi)$. Then
$\Sc^*_{(V\times V^*)(F)}(C(F))^{F^{\times}\times F^{\times}, \mu
\times 1 }=0.$
\end{lemma}
We will prove this lemma for Archimedean $F$, since the proof for
non-Archimedean is similar but simpler.

By Theorem \ref{NashFilt} this lemma follows from the following
one.

\begin{lemma}
For any $k \geq 0$ we have\\
(i) $\Sc^*(((V-0) \times 0) (F), Sym^k(CN_{(V-0) \times
0}^{V\times V^*}(F)))^{F^{\times}\times F^{\times}, \mu
\times 1 } = 0$.\\
(ii) $\Sc^*((0 \times (V^*-0)) (F), Sym^k(CN_{0 \times
(V^*-0)}^{V\times V^*}(F)))^{F^{\times}\times F^{\times}, \mu
\times 1 } = 0$.\\
(iii) $\Sc^*(0, Sym^k(CN_{0}^{V\times V^*}(F)))^{F^{\times}\times
F^{\times}, \mu
\times 1 } = 0$.\\
\end{lemma}

\begin{proof}$ $\\
(i) Cover $V-0$ by standard affine open sets $V_i:=\{x_i \neq
0\}$. It is enough to show that $\Sc^*((V_i \times 0) (F),
Sym^k(CN_{(V_i \times 0)(F)}^{V\times V^*}(F)))^{F^{\times}\times
F^{\times}, \mu \times 1 } = 0$.

Note that $V_i$ is isomorphic as an $F^{\times} \times
F^{\times}$-manifold to $F^{n-1} \times F^{\times}$ with the
action given by $(x,y)(v,\alpha) = (v,\frac{y}{x}\alpha)$. Note
also that the bundle $Sym^k(CN_{(V_i \times 0)(F)}^{V\times
V^*}(F))$ is a constant bundle with fiber $Sym^k(V)$.

Hence by Proposition \ref{Product} it is enough to show that
$\Sc^*(F^{\times}, Sym^k(V))^{F^{\times} \times F^{\times}, \mu
\times 1}=0$. Let $H := (F^{\times} \times F^{\times})_{1} = \{(t,
t) \in  F^{\times} \times F^{\times} \}$. Now by Frobenius descent
(Theorem \ref{Frob}) it is enough to show that
$(Sym^k(V^*(F))\otimes_{\R}\C) ^{H,\mu \times 1|_H}=0$. This is
clear since $(t,t)$ acts on $(Sym^k(V^*(F))$ by multiplication
by $t^{-2k}$.\\
(ii) is proven in the same way.\\
(iii) is equivalent to the statement $((Sym^k(V\times V^*)(F))
\otimes_{\R}\C) ^{F^{\times}\times F^{\times}, \mu \times 1 } =
0$. This is clear since $(t,1)$ acts on $Sym^k(V\times V^*)(F)$ by
multiplication by $t^{-k}$.
\end{proof}

\section{Proof of the geometric statement} \label{ProofGeo}
\begin{notation}
Denote $\check{S}'':=\{((A_1,v_1,\phi_1),(A_2,v_2,\phi_2)) \in
\check{S}'| A_1^{n-1}=0\}$.
\end{notation}
By Theorem \ref{non_disting_no_cois} and Example \ref{group_case}
there are no non-empty $X \times X$-weakly coisotropic
subvarieties of $\check{S}''$. Therefore it is enough to prove the
following Key proposition.
\begin{proposition}[Key proposition]
There are no non-empty $X \times X$-weakly coisotropic
subvarieties of $\check{S}'- \check{S}''.$
\end{proposition}
\begin{notation}
Let $A \in sl(V)$ be a nilpotent Jordan block. Denote $$R_A :=
(\check{S}'-\check{S}'')|_{\{A\} \times V \times V^*}.$$
\end{notation}
By Proposition \ref{GeoFrob} the Key proposition follows from the
following Key Lemma.
\begin{lemma}[Key Lemma]
There are no non-empty $V \times V^* \times V \times V^*$-weakly
coisotropic subvarieties of $R_A$.
\end{lemma}
\begin{proof}
Denote $Q_A = \bigcup_{i=1}^{n-1} (Ker A^i) \times
(Ker(A^*)^{n-i}).$ It is easy to see that $R_A \subset Q_A \times
Q_A$ and
$$Q_A \times Q_A = \bigcup_{i,j=0}^n (Ker A^i) \times (Ker (A^*)^{n-i}) \times (Ker A^j) \times (Ker (A^*)^{n-j}).$$
Denote $L_{ij} := (Ker A^i) \times (Ker (A^*)^{n-i}) \times (Ker
A^j) \times (Ker (A^*)^{n-j})$.

It is easy to see that any weakly coisotropic subvariety of $Q_A
\times Q_A$ is contained in $\bigcup_{i=1}^{n-1} L_{ii}.$ Hence it
is enough to show that for any $0<i<n$, we have $\dim R_A \cap
L_{ii} < 2n$.

Let $f \in {\mathcal O}(L_{ii})$ be the polynomial defined by
$$f(v_1,\phi_1,v_2,\phi_2):= (v_1)_{i} (\phi_2)_{i+1}-(v_2)_{i} (\phi_1)_{i+1},$$
where $( \cdot)_i$ means the i-th coordinate. It is enough to show
that $f(R_A \cap L_{ii}) = \{0\}$.

Let $(v_1,\phi_1,v_2,\phi_2) \in L_{ii}$. Let $M:= v_1 \otimes
\phi_2-v_2 \otimes \phi_1$. Clearly, $M$ is of the form
$$ M= \begin{pmatrix}
  &0_{i \times i} &* \\
  & 0_{(n-i) \times i} &0_{(n-i) \times (n-i)}
\end{pmatrix}.  $$
Note also that $M_{i,i+1}=f(v_1,\phi_1,v_2,\phi_2)$.

It is easy to see that any $B$ satisfying $[A,B]=M$ is upper
triangular. On the other hand, we know that there exists a
nilpotent $B$ satisfying $[A,B]=M$. Hence this $B$ is upper
nilpotent, which implies $M_{i,i+1}=0$ and hence
$f(v_1,\phi_1,v_2,\phi_2)=0$.

To sum up, we have shown that $f(R_A \cap L_{ii} = \{0\}$, hence
$dim(R_A \cap L_{ii})<2n$. Hence every coisotropic subvariety of
$R_A$ has dimension less than $2n$ and therefore is empty.
\end{proof}

\newpage
$ $
\newpage


\section*{Supplementary material}

משפט **ד**. כל התפלגות על (F⊗V)O שנשמרת על ידי הצמדה ב- O(V) נשמרת גם על ידי שחלוף.

כמו כן, לכל הצגה מותרת בלתי פריקה π של (F⊗V)O ולכל הצגה מותרת בלתי פריקה ρ של O(V) מתקיים

$$\dim \mathrm{Hom}_{O(V)}(\,\pi|_{O(V)},\,\rho)\le 1$$

עבור F לא ארכימדי משפט ד' הוכח ב-[AGRS] ועבור F לא ארכימדי משפט ד' הוכח ב-[SZ]. שתי ההוכחות השתמשו במשפט א'.

בהוכחות שלנו עבור F לא ארכימדי השתמשנו בכלים העוסקים עם התפלגויות אינוואריאנטיות שפותחו ב-

[BZ], [Ber], ו-[JR]. בשביל F ארכימדי פיתחנו מקבילים ארכימדיים של כלים אלה (ראה פרק 1).

המקרה הארכימדי של משפט א' הרבה יותר קשה ולכן הזדקקנו לכלי נוסף כדי לפתור אותו. הכלי המכריע היה התורה של D-מודולים. בפרט השתמשנו במשפט עמוק על התומך הסינגולרי של D-מודול. כדי שההוכחות שלנו יהיו אחידות לשדות ארכימדיים ולא ארכימדיים השתמשנו במאמר [Aiz] שפותח כלי מקביל למקרה הלא ארכימדי.

# סיכום

יהי F שדה מקומי מאופיין אפס, לדוגמא שדה המספרים הממשיים או שדה המספרים ה-p-אדיים. נתבונן בשיכון הטבעי GL(n,F)<GL(n+1,F) . נתבונן בפעולת ההצמדה של GL(n,F) על GL(n+1,F). התוצאה העיקרית של תזה זו היא המשפט הבא.

משפט **א**. כל התפלגות על GL(n+1,F), אשר נשמרת על ידי פעולת ההצמדה של GL(n,F) נשמרת גם על ידי שחלוף.

הוכחנו זאת ב- [AGRS] עבור F לא ארכימדי וב-[AG08e] עבור F ארכימדי. עבור F ארכימדי, המשפט הוכח גם ב-[SZ] בו זמנית, באופן בלתי תלוי ובדרך אחרת. ההוכחה שאנו מביאים בתזה זו היא שילוב של שלושת ההוכחות הנ"ל ועובדת בצורה אחידה עבור שדות ארכימדיים ולא ארכימדיים.

למשפט חשיבות בתורת ההצגות משום שהוא גורר את משפט חוסר ריבוי הבא

משפט **ב**. תהי π הצגה מותרת בלתי פריקה של GL(n+1,F) ו ρ הצגה מותרת בלתי פריקה של GL(n,F). אזי

$$\dim \mathrm{Hom}_{GL(n,F)}(\pi|_{GL(n,F)}, \rho) \leq 1$$

מסקנה נוספת של משפט א' היא

משפט **ג**. נסמן ב P(n,F) < GL(n,F) את התת-חבורה המורכבת ממטריצות שהשורה האחרונה שלהן היא (0,...,0,1). אזי כל התפלגות על GL(n,F) שנשמרת על ידי הצמדה ב- P(n,F) נשמרת על ידי הצמדה ב- GL(n,F).

משפט ג' גורר את השארת Kirillov שטוענת כי כל הצגה אוניטרית בלתי פריקה של GL(n,F) נשארת בלתי פריקה כהצגה אוניטרית לאחר צמצום ל- P(n,F). בדרך זו השארת Kirillov הוכחה במקור, ב-[Ber] למקרה הלא ארכימדי וב-[Bar] למקרה הארכימדי.

משפטים דומים למשפטים א' וב' מתקיימים גם עבור חבורות אורתוגונליות. יהי V מרחב וקטורי ממימד סופי מעל F ותהי q תבנית ריבועית לא מנוונת על V. נרחיב את q בדרך הטבעית ל- V⊕F ונתבונן בשיכון O(V) <O(V⊕F).

לזכר הסבים שלי

נאום גורביץ ודוד קגן

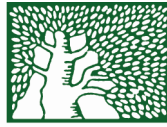

WEIZMANN INSTITUTE OF SCIENCE

Thesis for the degree
**Doctor of Philosophy**

Submitted to the Scientific Council of the
Weizmann Institute of Science
Rehovot, Israel

עבודת גמר (תזה) לתואר
**דוקטור לפילוסופיה**

מוגשת למועצה המדעית של
מכון ויצמן למדע
רחובות, ישראל

במתכונת "רגילה"
In a "Regular" Format

By
**Dmitry Gourevitch**

מאת
דמיטרי גורביץ

משפטי חוסר ריבוי
והתפלגויות אינוואריאנטיות
Multiplicity One Theorems
and Invariant Distributions

Advisors:
Prof. Stephen Gelbart  and
Prof. Joseph Bernstein

מנחים:
פרופ' סטפן גלברט
ופרופ' יוסף ברנשטיין

May 2009

סיוון תשס"ט


\end{document}